\documentclass{article}
\pdfoutput=1
\usepackage{arxiv}
\usepackage[utf8x]{inputenc} 
\usepackage[T1]{fontenc}    
\usepackage{caption}
\usepackage{subcaption}
\usepackage{graphicx}
\usepackage{color}
\usepackage[section]{placeins}
\usepackage[normalem]{ulem}
\usepackage{amsmath,amsbsy,amssymb}
\usepackage{hyperref}
\usepackage{tikz}
\usepackage[title]{appendix}

\title{An explicit semi-Lagrangian, spectral method for solution of Lagrangian transport equations in Eulerian-Lagrangian formulations}

\author{
Hareshram Natarajan \\
Computational Science Research Center \\
San Diego State University\\
San Diego\\
\AND
Gustaaf B. Jacobs\\
Department of Aerospace Engineering\\
San Diego State University\\
San Diego\\
}

\begin{document}
\maketitle

\begin{abstract}
An explicit high order semi-Lagrangian method
is developed for solving  Lagrangian transport equations in Eulerian-Lagrangian
formulations. To ensure a semi-Lagrangian approximation 
that is consistent with an explicit Eulerian, discontinuous spectral element method (DSEM) discretization used for the Eulerian formulation, Lagrangian particles are seeded at Gauss quadrature
collocation nodes within an  element. The particles are integrated explicitly in time
to obtain an advected  polynomial solution at the advected
Gauss quadrature locations. This approximation is mapped back 
in a semi-Lagrangian fashion to the
Gauss quadrature points through a least squares fit using constraints for
element boundary values and optional constraints for mass and energy
preservation. An explicit time integration is used
for the semi-Lagrangian approximation that is consistent
with the grid based DSEM solver, which  ensures that particles seeded at the Gauss quadrature points do not leave the
element's bounds. The method is hence local and parallel and
facilitates the solution of the Lagrangian formulation without
the grid complexity, and parallelization challenges of a particle solver in particle-mesh methods. 
Numerical tests with one and two dimensional advection equation are carried out.
The method converges exponentially. The use of mass and energy
constraints can improve accuracy depending on the accuracy of the time integration.
\end{abstract}


\section{Introduction}
\label{intro}
A range of problems can be described by 
coupled Eulerian-Lagrangian (EL) formulations 
where a hyperbolic Eulerian governing system of equations such as the Euler equations
for gas dynamics or the Maxwell's equations for electromagnetics is coupled
with a Lagrangian governing system for particles such as solid particulates or
charged particles. Interestingly, either one of these Eulerian and
Lagrangian systems can often be recast into an
equivalent Lagrangian or Eulerian formulation, respectively.
The kinetic, Lagrangian behavior of solid particles or charged particles,
for example, can be cast into an Eulerian formulation through
a method of moments or distributions leading to Boltzmann
Vlasov and/or Fokker-Planck type equations. Vice-versa the Eulerian
gas dynamics have an equivalent Lagrangian form along the path
of a fluid particle.

In the approximation of the EL formulations
a similar duality exists.  While it is common that Eulerian equations are approximated
with grid based methods such as finite volume or finite difference methods,
Lagrangian formulations are usually solved with particle methods.
Particle-mesh methods which couple the two approaches are predominant
for the solution of EL systems. 
However, one could also consider Lagrangian
methods, such as point vortex methods \cite{Chorin} to solve
the Eulerian systems, or vice-versa grid-based methods to solve
the Lagrangian system using its equivalent Eulerian form \cite{ShotorbanJacobs}.
For the interested reader, we refer to text books such as Birdsall et al. \cite{Birdsall}
and Crowe et al. \cite{Crowe} for a comprehensive review of classic formulations
and methods.

High-order accurate schemes like discontinuous spectral element methods (DSEM) \cite{KoprivaJacobs,Kopriva}
are a particularly good choice to solve hyperbolic, Eulerian
equations, because they propagate waves over long distances
and capture small scales better than low-order methods. They do so with dramatically
reduced resource needs, e.g. being accurate on coarser grids or having high efficiency at high accuracy.
Because of their local nature, DSEM Navier-Stokes solvers have been shown extensively to
obtain high accuracy and convergence using unstructured grids on complex geometries
 \cite{HW08, Kopriva09}.
Moreover, both in theory \cite{HW08}
and in testing through benchmarks (e.g. \cite{JKM04,Kopriva,HOworkshop}),
DSEMs have been shown to have superior computational efficiency and parallelism
for computation of smooth flows as compared to more traditional discretization methods.

Particle methods for Lagrangian formulations
are inherently conservative. Lagrangian methods typically have 
a time step criterion that is less restrictive for stability than numerical
methods based on the equivalent Eulerian formulations \cite{lagpart}. A downside of Lagrangian
tracers is that they are dynamically moving in a given
computational domain and are hence less tractable than
static Eulerian grid points \cite{atm1}. 

To preserve the favorable properties of the DSEM method
in an EL formulation, the coupling
with the grid based solver and the Lagrangian, particle solver
must be consistently high-order. In Refs. \cite{KoprivaJacobs,JD09,JH09,JH06,SS11,SJ16,SJD14}, 
a consistent high-order interpolation of the Eulerian
solution to the particle and a consistent high-order distribution
of the particle influence on to the grid was developed.
The tracking of the connectivity between the dynamic particle mesh and the static grid constitutes
a major challenge. It increases code complexity and is detrimental to parallel efficiency.
For general complex grids, ensuring a consistent coupling
and formal convergence is a topic of ongoing research.  

A semi-Lagrangian (SL) method has the potential to take advantage of the 
best of the Eulerian and Lagrangian worlds. In SL methods, the
solution can be traced either forward or backward in 
time along characteristic lines. The traced solutions is then remapped
onto an underlying static grid  after each time step \cite{atm4,remap1,atm1}. 
While for simple implementations SL is typically non-conservative, 
there have been several conservative
remapping strategies proposed. In \cite{remap2,remap3} conservative remapping
of quadrilateral elements is achieved by the application of Gauss-Green’s
theorem which converts area-integrals into line-integrals. \cite{cslam}
extended this work to cubed-sphere grids. The use of straight lines to
approximate control volume limits these methods to second order.

In the semi-Lagrangian flux form,
only boundary fluxes are treated in a semi-Lagrangian manner. 
Starting with the Eulerian form of the transport equation, 
a semi-Lagrangian time integration 
evaluates the fluxes at element boundaries\cite{flux1,atm2,atm3}. 
This method implicitly guarantees
conservation. In order to extend the scheme to two and three dimensions, a
dimensional splitting technique is used which limits the method to second
order.

Semi-Lagrangian have been developed for continuous and discontinuous spectral element methods
\cite{dgsl2,dgsl2b,dgsl4} to solve Vlasov or Boltzmann type equations.
The first semi-Lagrangian spectral
element method was developed in \cite{dgsl2} based on Lagrange-Galerkin formulation. The method is shown to be stable
and high order accurate, but is not strictly conservative. This work is extended to
shallow water application in \cite{dgsl2b}. Flux form semi-Lagrangian method in
discontinuous Galerkin(SLdG) framework was introduced in \cite{flux2}. The
method is shown to be stable, conservative and show results of third order
accuracy. Similar methods with higher order accuracy were developed in
\cite{dgsl4,dgsl3}. Recently, a conservative SLdG method which tracks a
Lagrangian control volume has been developed in \cite{vp2}. The method is shown
to be stable with fifth order accuracy in space for uniform velocity problems.
\cite{dgsl1} shows a similar method that permits non uniform velocities. Most
of the current SLdG methods require backward tracking of element interfaces to
calculate the fluxes. This makes parallel scalability difficult when the
geometry is complex. 

To the best of our knowledge, all DSEM formulations developed
so far rely on the flux based form of the conservation law.  
No DSEM based semi-Lagrangian method that starts from the Lagrangian form 
of the equations in an EL formulation
has been developed.

In this paper, we develop an explicit high order semi-Lagrangian method 
that is consistent with the
DSEM framework. The SL method solves the Lagrangian form of the transport equations and is intended to be  used instead of the particle solver in particle-mesh methods.
Because Eulerian DSEM is local
and highly parallelizable, we 
present a consistent SL methods that is local as well.
The size of linearized systems for an Euler solver coupled
with a particle method is immense and rather complex and challenges
implicit methods and effective iterative solver. 
For coupled EL systems, an explicit time-integrator is a natural choice
and most commonly used. 
We hence focus on explicit time-integrators and ensure
that the SL method is explicit also.
By seeding Lagrangian particles 
on Gauss collocation nodes within a spectral element,
we  are able to develop a local and explicit SL method
with spectral resolution. The particles seeded at the Gauss quadrature points are 
integrated forward in time along their characteristic path. 
The time step is restricted such that the particles do not
cross the element boundaries and hence do not require parallel attention other
than at the elements interfaces when the advected solution is patched at
interfaces. The final step is to remap the
solution at the collocation nodes at the new time using a least square
fit with boundary constraints. The method is not inherently conservative. Additional constraint equations to ensure on local mass and energy are implemented and tested for the cases of a one-dimensional and two-dimensional constant advection velocity as well as a variable advection velocity.  

In the next section, we present discuss Lagrangian transport equations in Eulerian-Lagrangian formulations, which is  followed by an overview of the DSEM method. Then the DSEM-SL algorithm is presented. Numerical results for one and two dimensional tests are given next. Conclusions and future steps are reserved for the final section.

\section{Governing equations}
\label{sec:1}
\subsection*{General Eulerian-Lagrangian formulation}

We focus on Lagrangian transport equations in general
Eulerian-Lagrangian formulations. In the Eulerian frame, 
we consider hyperbolic  conversation laws of the following form:
\begin{equation} 
\frac{\partial {\mathbf{Q}}}{\partial t} + {\nabla}\cdot {\mathbf{F}} =\mathbf{S}=
	\sum_i^{Np} \mathbf{F^p_i} \left(\mathbf{Q}(\mathbf{x^p_i}),\mathbf{Q^p_i}\right)  K(\mathbf{x},\mathbf{x^p_i})
\label{eq:Eulerian}
\end{equation}
where $\mathbf{Q}$ and $\mathbf{F}$ are the solution and
flux vector, respectively. Through the source term, $\mathbf{S}$, these
equations are coupled to a Lagrangian formulation of the form
\begin{eqnarray}
	{d \mathbf{x^{p}_i} \over d t} &=& \mathbf{v^p_i} \label{eq:Lagrangian_xp} \ \ \ \ i=1...Np \\
	{d \mathbf{Q^{p}_i} \over d t} &=& \mathbf{F^p_i}(\mathbf{Q}(\mathbf{x^p_i}),\mathbf{Q^p_i})
\label{eq:Lagrangian}
\end{eqnarray}
where $\mathbf{x^p_i}$ identifies the $i^{th}$ particle's location along its
Lagrangian path and $\mathbf{v^p_i}$ is its velocity. The particle's solution
$\mathbf{Q^p_i}$ is advected along the particle's path and
is forced according to $\mathbf{F^p_i}(\mathbf{Q}(\mathbf{x^p_i}),\mathbf{Q^p_i})$.
The forcing is a function of both the Eulerian and the Lagrangian solution
and is responsible for the coupling from the Eulerian to the Lagrangian frame.
The source term in (\ref{eq:Eulerian}), which inversely couples
the Lagrangian to the Eulerian equation, is also a function of both solutions
and is the summation of each particle's forcing distributed over space
according to a distributing function, $K(\mathbf{x},\mathbf{x^p_i})$.

\subsection*{Lagrangian Transport Equations; Reduced Form}

To discuss the development of the semi-Lagrangian method,
we focus on the Lagrangian transport equations and therefor consider a reduced passive Lagrangian formulation that does not couple back to the Eulerian formulation.
Moreover, we take the particle's position according
(\ref{eq:Lagrangian_xp}) to be determined
by a prescribed velocity field, $\mathbf {u}$, that is
obtained from the Eulerian solution, i.e.
the velocity, $v(\mathbf{Q})$, is a function of $\mathbf{Q}$. It is
quite common for this type of tracer particle to be used in Eulerian-Lagrangian formulations (\cite{Jaberi}).
In general, the scalar advection equation describing transport phenomena in the conservation form is given by, 
\begin{equation} \label{eq:2.3}
\frac{\partial \phi}{\partial t} + \nabla.(u \phi) =0, 
\end{equation}
where $\phi$ is the variable which is  transported. This equation can be written in the non-conservation form
as follows,
\begin{equation} \label{eq:transport_noncons}
\frac{\partial \phi}{\partial t} + \mathbf{u}.\nabla \phi =  -\phi(\nabla.u). 
\end{equation}
In this non-conservative form, the right hand side of the equation is usually zero in hyperbolic
conservation laws as the mass conservation (in \ref{eq:transport_noncons}, the divergence of the velocity) factors out. 

In Lagrangian form, along characteristic curves, this can be equivalently formulated as a combination of a kinematic
equation along the characteristic path as
\begin{equation} \label{eq:2.5}
\frac{d \mathbf{x^p}(t)}{d t} = \mathbf{u}(\mathbf{x^p}(t),t),
\end{equation}  
with an equation for the advected solution as follows, 
\begin{equation} \label{eq:2.4}
\frac{D \phi}{D t} = -\phi(\nabla.u). 
\end{equation}
Here, we focus on (\ref{eq:2.5}) and (\ref{eq:2.4})
to introduce the semi-Lagrangian method for solution
of Lagrangian equations. These equations are simplified, but
equivalent forms of (\ref{eq:Lagrangian})
with $\mathbf{v}$ =$\mathbf{u}$, 
$\mathbf{Q_i^p}$= $\rho \phi$ and a  zero forcing 
$\mathbf{F^p_i(\mathbf{Q}(\mathbf{x^p_i}),\mathbf{Q^p_i})}$=0.


\section{Discontinuous Spectral Element Method (DSEM) }
\label{sec:2}
We develop a semi-Lagrangian method for Eulerian-Lagrangian formulations 
that is consistent with the staggered grid DSEM  approximation 
as first introduced by Kopriva \cite{Kopriva}
for approximation of the  Eulerian equations.
In this version of DSEM the  solution variable is collocated
at Gauss quadrature nodes and the fluxes on Lobatto quadrature nodes.
The collocation at Gauss quadrature nodes turns out to be specifically beneficial
to preserve the local nature of DSEM in the SL formulation
as we shall see below.  To introduce notation and to set the stage,  we briefly
summarize essential aspects of  the staggered grid DSEM method.  For a detailed
description, we refer to \cite{Kopriva09,JKM04}.

In DSEM, the physical domain $\Omega$ is divided
into $K$ non-overlapping elements, $\Omega = \cup_{k=1}^K \Omega_k $. 
In the context of DSEM, elements are often referred to as subdomains, a nomenclature
that we follow in this paper. Each physical
subdomain is then mapped onto a unit computational cube using iso-parametric
transformation \cite{KoprivaJacobs}. This transforms the governing Eulerian equation, (\ref{eq:Eulerian}), to
\begin{equation} 
\frac{\partial \tilde{\mathbf{Q}}}{\partial t} + \tilde{\nabla}\cdot\tilde{F} =0,
\label{eq:3.1}
\end{equation}
where, $\tilde{\mathbf{Q}}= |\overline{\overline{J}}| \mathbf{Q}$, $\tilde{\nabla}\cdot\tilde{F} = \frac{\partial \tilde{f}}{\partial \xi} + \frac{\partial \tilde{g}}{\partial \eta} + \frac{\partial\tilde{h}}{\partial \zeta}$. $|\overline{\overline{J}}|$ is the determinant of the transformation from the physical to the computational domain. 

The solution and flux collocation points are chosen according to Chebyshev Gauss and Lobatto quadrature points,
which along tensorial grid lines,  $0\leq \xi \leq 1$, are given by,
\begin{equation} \label{eq:gaussquad}
\xi_{i+1/2} = \frac{1}{2}\left[ 1- \cos\left(\frac{i+1/2}{N+1}\right) \pi \right] \qquad i=0,1,...,N-1,
\end{equation}
and
\begin{equation} 
\xi_i = \frac{1}{2}\left[1-\cos\left(\frac{i \pi}{N}\right)\right] \qquad i=0,1,...,N,
\label{eq:3.3}
\end{equation} 
respectively.
Here, we have used the integer subscript, $i$, to identify Lobatto
points and $i+1/2$ to identify Gauss points that are located
in between two Lobatto points $i$ and $i+1$.
In three dimensions, the solution interpolant $\tilde{\mathbf{Q}}$ is then
\begin{equation} 
\tilde{\mathbf{Q}}(\xi, \eta, \zeta) = \sum \limits_{i=0}^{N-1} \sum \limits_{j=0}^{N-1} \sum \limits_{k=0}^{N-1} {\tilde{\mathbf{Q}}}_{i+1/2,j+1/2,k+1/2} h_{i+1/2}(\xi) h_{j+1/2}(\eta) h_{k+1/2}(\zeta),
\label{eq:3.4}
\end{equation}
where $h_{i+1/2}(\xi)$ is the Lagrange interpolation polynomial of degree N-1 defined on the Gauss 
quadrature points $\xi_{m+1/2}$ 
and
\begin{equation} 
h_{i+1/2}(\xi) = \prod_{\substack{m=0 \\ m \neq p}}^{N-1} \frac{\xi - \xi_{m+1/2}}{\xi_{i+1/2} - \xi_{m+1/2}}, \qquad i=0,1,...,N-1,
\label{eq:3.5}
\end{equation}
is the Lagrangian polynomial of degree $N$-1.
The fluxes, $\tilde{F}$, are collocated  similarly on the Lobatto points. Through interpolation 
between the Gauss grid and the Lobatto grid, the fluxes can be determined as a function
of the solution, $\tilde{Q}$.  Through an approximated Riemann solver,
an interface flux is determined from interface solutions on neighbouring subdomains.
The derivatives of the  fluxes, $\tilde{\nabla}\cdot\tilde{F}$, are  determined at the Gauss points.
Then, it remains to update the Gauss solution  in time. We typically use an
explicit integrator such as a  standard fourth order explicit Runge-Kutta time stepping method.

\section{Semi-Lagrangian scheme for Lagrangian transport equations}
\label{sec:3}
\begin{figure}
	\centering 
	\mbox{
		\includegraphics[width=0.75\textwidth]{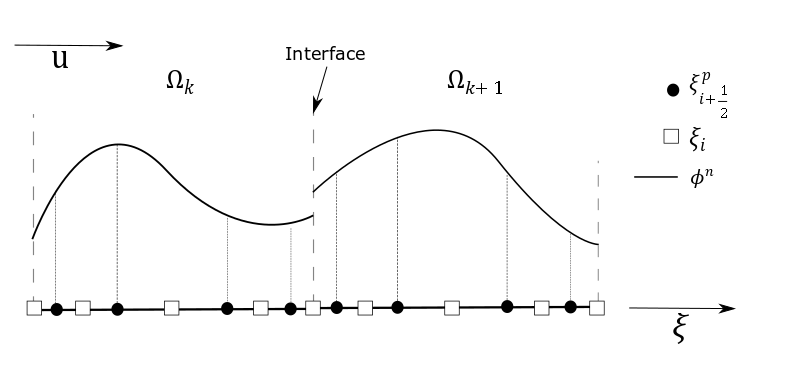}
	}
	\makebox[0.65\textwidth][c]{(a) Initialization at Gauss points}
	\mbox{
		\includegraphics[width=0.75\textwidth]{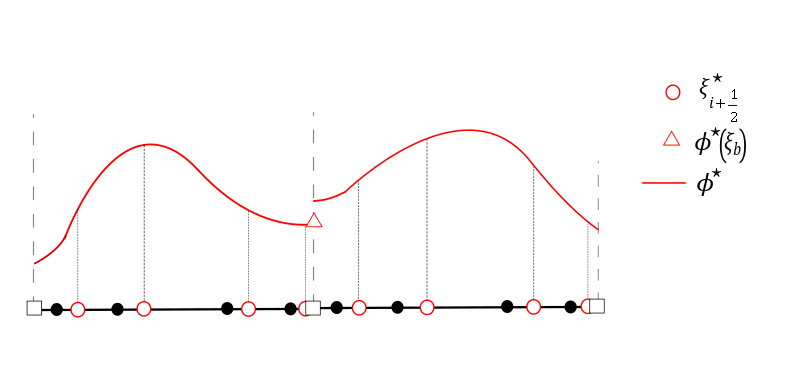}
	}
	\makebox[0.65\textwidth][c]{(b) Particle advection}
	\mbox{
    	\includegraphics[width=0.75\textwidth]{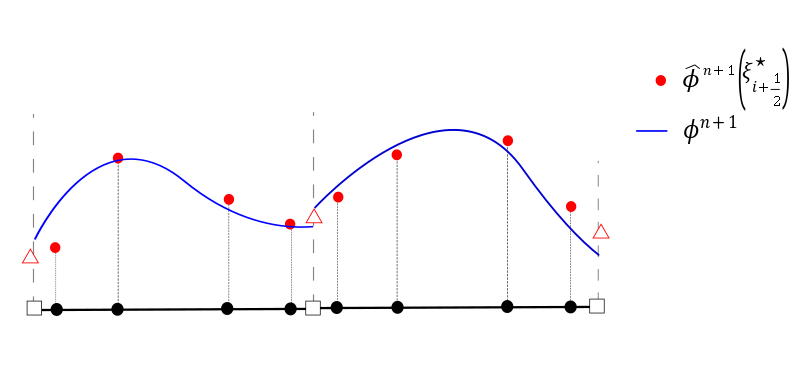}
	}
	\makebox[0.65\textwidth][c]{(c) Remapping to Gauss points}
	\caption{Schematic of the DSEM SL method: (a) The tracer particles are initialized on each subdomain at the Chebyshev Gauss points, $\xi_{i+1/2}^p$; (b)  the particles are traced forward in time to obtain the advected location, $\xi^\star_{i+1/2}$, the advected solution, $\phi^\star$, and the the interface solution, $\phi^\star(\xi_b)$; (c) the solution is remapped to obtain the intermediate solution at original tracer locations, $\hat{\phi}^{n+1}(\xi_{i+1/2}^p)$. Using a least squares fit with constraint for the boundary condition, and for mass and/or energy conservation, we obtain the solution $\phi^{n+1}$.}
	\label{fig:0_1}
\end{figure}

The semi-Lagrangian algorithm that we propose requires
a number of steps that are shown
in a schematic in Figure\ \ref{fig:0_1}.
They include initialization
of particles at the Gauss points
(Fig.\ \ref{fig:0_1}a), 
advection of the particles through a time integration (Fig.\ \ref{fig:0_1}b) and a remapping 
of the advected  particle solution back to the Gauss points (Fig.\ \ref{fig:0_1}c). We describe
each step in detail below for a one-dimensional approximation.
The multidimensional algorithm mostly extends naturally through
a tensorial grid, except for the remapping which we discuss
separately.

\subsection{Solution initialization}
To be consistent with an Eulerian solver, 
we initialize $N$ particles within a subdomain, $k$, at  time $t^n$ at the 
Gauss quadrature collocation points in (\ref{eq:gaussquad}). 
The solution, $\phi^n$, is then approximated by a Lagrange interpolant as follows,
\begin{equation} \label{eq:4.2}
	\phi^n(\xi) = \sum \limits_{i=0}^{N-1}{\phi}^n(\xi^p_{i+1/2}) h_{i+1/2}(\xi), 
\end{equation} 
where $h_{i+1/2}(\xi)$ are the Lagrange polynomials of degree $N-1$ defined on the Chebyshev Gauss points $\xi_{i+1/2}^p$ according to \label{eq:3.2}
\begin{equation} \label{eq:4.3}
h_{i+1/2}(\xi) = \prod_{\substack{i=0 \\ i \neq j}}^{N-1} \frac{\xi - \xi^p_{i+1/2}}{\xi^p_{j+1/2} - \xi^p_{i+1/2}}, \qquad j=0,1,...,N-1.
\end{equation}
The solution $\phi$ is initialized at $t^0$=0 as ${\phi}^0(\xi^p_{i+1/2})$.

\subsection{Forward time integration}
The particle variable, $\phi^n$, is advected in time
along its characteristic path according to (\ref{eq:2.5}).
We update each particle's location, i.e.
\begin{equation} \label{eq:4.4}
	\frac{d \xi^p_{i+1/2} }{dt} =  u(\xi^p_{i+1/2})  \qquad i=0,1,...,N-1,
\end{equation}
with an explicit time integration,  so
\begin{equation} \label{eq:4.4a}
	\xi^\star_{i+1/2} =  	\xi^p_{i+1/2} + \Delta t u(\xi^p_{i+1/2})  \qquad i=0,1,...,N-1,
\end{equation} 
where $\xi^\star$ is the advected particle location at $t^{n+1}$.
In this paper, we take $u(\xi_{i+1/2})$ prescribed, but in a general Eulerian-Lagrangian
formulations it is a velocity interpolated from the Eulerian
grid based solver. Because the particles are seeded
at the Gauss quadrature points, the velocity at those
location is conveniently and directly available (without interpolation) from the DSEM solver.

An explicit time integration  is common for Eulerian-Lagrangian methods \cite{JH06,JH09, JD09}
In principle, any explicit integration scheme can be implemented for the DSEM
solver and consistently for the semi-Lagrangian solver also.
For simplicity of notation and for the sake of the 
concise presentation of the algorithm, we describe the semi-Lagrangian
method using a first-order Euler scheme in (\ref{eq:4.4a}).
We have used several time integration schemes of higher order in the tests discussed 
in the sections below. 

While there is no formal stability criterion for the temporal update
of the linear characteristic equation, in order  
to prevent an advected particle from leaving a subdomain
we restrict the  time step nevertheless as follows
\begin{equation} \label{eq:1d_delt_condition}
	\Delta t \leq \frac{\Delta \xi_{\text{min}}}{|u|_{max}}.
\end{equation}
Here, $\Delta \xi_{min}$ is the minimum grid
spacing between two particles, i.e. Chebyshev quadrature points, at the edges of the subdomain
and $|u|_{max}$ is the maximum advection speed. This time step limit is typically less restrictive than the explicit CFL time step for the Eulerian DSEM solver, because
the latter time step is restricted by the same minimum grid  spacing, $\xi_{min}$, but by
the maximum of multiple characteristic velocities of the hyperbolic system.
These velocities usually have a greater magnitude than $|u|_{max}$ and
so the explicit time step for the semi-Lagrangian method is not-restrictive for a coupled EL formulation. 

The benefit of preventing a particle from leaving the subdomain is 
two-fold. Firstly, information from the neighboring subdomains are not involved in the advection
step and the remapping will require the solutions to be
connected along the interface between subdomains only, i.e.
the subdomains are non-overlapping.  This yields  a local method; again, consistent with DSEM.
Secondly, with the advected particle locations within a subdomain deviating only marginally
from the quadrature points, we find that the advected polynomial is not ill-posed.


The solution after advection is denoted by $\phi^\star(\xi)$ and is given as
\begin{equation}
	\phi^\star(\xi) = \sum \limits_{i=0}^{N-1}{\phi^\star}(\xi^{\star}_{i+1/2}) h^{\star}_{i+1/2}(\xi), 
\label{eq:4.6}
\end{equation} 
where $h^{\star}_{i+1/2}(\xi)$ are the Lagrange  polynomials of degree $N-1$ 
defined on the advected points $\xi_{i+1/2}^{\star}$,
\begin{equation} 
h^{\star}_{i}(\xi) = \prod_{\substack{i=0 \\ i \neq j}}^{N-1} \frac{\xi - \xi^{\star}_{i+1/2}}{\xi^{\star}_{j+1/2} - \xi^{\star}_{i+1/2}}, \qquad j=0,1,...,N-1.
\label{eq:4.7}
\end{equation}
The advected polynomial's nodal solution values, ${\phi^\star}(\xi_{i+1/2}^{\star})$ are obtained by integrating $\phi^n(\xi_{i+1/2}^n)$ along the characteristic curve according to (\ref{eq:2.4}), 
\begin{equation} \label{eq:4.8}
	{\phi^\star}(\xi^{\star}_{i+1/2}) = {\phi}^n(\xi^p_{i+1/2}) + \Delta t \left( -{\phi}^n(\xi^p_{i+1/2}) \left(\frac{\partial u}{\partial \xi}\right)^n_{i+1/2}\right).
\end{equation}
Here, we take $\left(\frac{\partial u}{\partial \xi}\right)^n_{i+1/2}$ prescribed from a hypothetical Eulerian grid based solver.
\subsection{Remapping}
Finally, in a remapping stage the advected polynomial is projected back onto the Gauss-Chebyshev quadrature nodes through interpolation
%
as follows,
\begin{equation} \label{eq:4.9}
	\phi^\star(\xi^p_{i+1/2}) = \sum \limits_{j=0}^{N-1}{\phi^\star}(\xi^{\star}_{j+1/2}) h^\star_{j+1/2}(\xi^p_{i+1/2}), 
\end{equation}
 providing an estimate for the solution at time step, $n+1$, as follows
\begin{equation}
	\widehat{\phi}^{n+1}(\xi_{j+1/2}^p) = \phi^\star(\xi^p)  \qquad j=0,1,...,N-1.
\label{eq:remap_advect}
\end{equation}
Here, we use the \textit{hat} symbol to denote the intermediate solution at $t^{n+1}$.
To account for connectivity between elements, boundary conditions
and conservation of mass, 
we constrain this intermediate solution. We then use a least-squares
projection to obtain a corrected solution at $t^{n+1}$.

\begin{figure} 
	 \centering
		\includegraphics[width=0.75\textwidth]{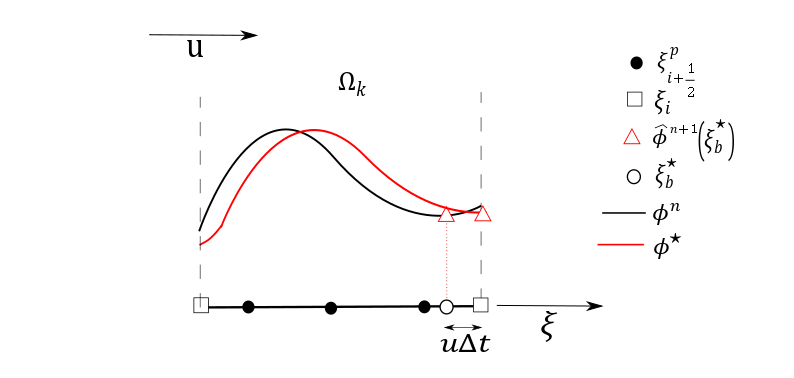}
	\caption{Schematic of backward tracing procedure to determine interface values of the solution: to determine the interface value $\hat{\phi}^{n+1}(\xi_b)$, a characteristic is traced backward in time that  passes through $\xi_b$.The origin of the characteristic is located using $\xi^\star=\xi_b -u\Delta t$. The solution is then interpolated at $\xi^\star_b$, using the known solution at the quadrature nodes at $t^n$, $\phi^n(\xi_{i+1/2})$ .}
	\label{fig:0_2}
\end{figure}

Boundary conditions and interface constraints can be applied either using backtracking 
on the characteristic line or by interpolation. 
To find the interface value $\hat{\phi}^{n+1}(\xi_b)$, we trace a characteristic passing through $\xi_b$ backward in time as is illustrated in
Figure\ \ref{fig:0_2}. We locate the origin of the characteristic using $\xi^\star=\xi_b -u\Delta t$. The solution is then interpolated at $\xi^\star_b$, using the known solution at the quadrature nodes at $t^n$, $\phi^n(\xi_{i+1/2})$. \cite{backtrack}.

\begin{eqnarray}
	\hat{\phi}^{n+1}_b & = &  \sum  \limits_{j=0}^{N-1}{\phi}^n(\xi^p_j) h_{j}(\xi_{b}-u\Delta t)  \qquad u>0 \qquad b=1,2. 
\label{eq:4.11}
\end{eqnarray}
where the subscript $b$=1,2 identifies an interface value or coordinate 
at the left or right end of the subdomain, respectively.

Alternatively, we can determine the boundary values using polynomial interpolation 
according to (\ref{eq:4.6}),
\begin{eqnarray} \label{eq:remap_bnd}
	\hat{\phi}^{n+1}_b = \sum \limits_{j=0}^{N-1}{\phi^\star}(\xi^{\star}_{j+1/2}) h^\star_{j+1/2}(\xi_b)  \qquad b=1,2 
\end{eqnarray}
By upwinding, a unique interface value  is determined
from the interfaces values of two neighbouring subdomains
\begin{eqnarray} \label{eq:4.12a}
	\hat{\phi}^{n+1}_{b\star} =  f \left(\hat{\phi}^{n+1}_{b=1} \bigg\vert_{\Omega_k}, \hat{\phi}^{n+1}_{b=2} \bigg\vert_{\Omega_{k-1}} \right)
\end{eqnarray}
Boundary conditions are implemented in the same way as interface
condition by using a specified ghost solution at computational domains boundaries.

In the tests conducted below, we compare the 
backtracing and interpolation approaches and find that they produce similar results.

\subsection*{Mass and Energy Conservation}
\subsubsection*{Local mass conservation}
With the above remapping, mass is not exactly conserved. An additional constraint for local mass conservation can be imposed as follows. Using
the equivalent Eulerian formulation in (\ref{eq:2.3}),
an expression for the mass in a subdomain is given by
\begin{equation} \label{eq:4.14}
	\frac{\partial}{\partial t} \int_{\Omega_k} \phi dx  = - u \left.  \phi \right|^{\xi_{b=2}}_{\xi_{b=1}}.
\end{equation} 
Here, $M$ is the mass, $M =\int_{\Omega_k} \phi dx$  and  
$\mathcal{F}_x$ is the mass flux, 
$\mathcal{F}_x =u \left. \phi \right|^{\xi_{b=2}}_{\xi_{b=1}}$.
The integral is numerically evaluated using Gauss quadrature,
\begin{equation} \label{eq:quad}
	M^n = \sum \limits_{j=0}^{N-1} w_{j+1/2} {\phi}^{n}(\xi_{j+1/2}) \qquad j =0,1,...,N-1,
\end{equation} 
where $w_{j+1/2}$ are the quadrature weights. 
We find the mass at $M^{n+1}$ through time integration, 
\begin{equation} \label{eq:remap_mass}
M^{n+1} = M^n - \mathcal{F}_m^n \Delta t,
\end{equation}
with the flux $\mathcal{F}_m^n = u \phi^n_{b=1} -  u \phi^n_{b=2} $
The boundary values, $\phi^n_b$, at $t^n$ are determined with interpolation
from the polynomial defined on the Gauss qaudrature nodes.

\subsubsection*{Local energy conservation}
Just like mass, energy is not exactly conserved. An additional constraint for local energy conservation with a constant velocity advection can be formulated as follows,
\begin{equation}
\frac{\partial}{\partial t} E = -u \left. \phi^2 \right|^{x_\text{b=2}}_{x_\text{b=1}},
\end{equation}
where $E=\int_{\Omega_k} \phi^2 dx$.
The integral is evaluated using Gauss quadrature. 
Because the energy is non-linear (quadratic
in $\phi$), we have to iterate to solve the constrained system
that will be discussed below.
For that purpose we introduce the $k$ counter as follows
\begin{equation}
	E^{n} =  \sum \limits_{j=0}^{N-1} w_{j+1/2} \widehat{\phi}^{(k-1)}(\xi_{j+1/2}) \widehat{\phi}^{(k)}(\xi_{j+1/2}) \qquad j=0,1,...,N-1,
\end{equation}
The iteration is terminated when $|\widehat{\phi}^{(k-1)}-\widehat{\phi}^{(k)}| < 10^{-8}$.
For a first order time integration and a constant $u$, the energy
at time $t^{n+1}$ is determined as,
\begin{equation} \label{eq:remap_energy}
	E^{n+1} = E^n -  \Delta t \mathcal{F}_e^n, 
\end{equation}
where $\mathcal{F}_e^n=u \left. \left(\phi^{n}\right)^2 \right|^{x_\text{b=2}}_{x_\text{b=1}}$.


\subsubsection*{Least-squares solution of the  constrained system}
To project the interpolated polynomial, $\hat{\phi^{n+1}}$, combined with the constraints
onto the Gauss-Chebyshev quadrature using a least-squares method,
an overdetermined system of  equations has to be solved. The number of equations depends on the number of constraint equations included. For example, if the boundary constraint and the local mass and energy constraint equations are used, Equations
(\ref{eq:remap_advect}), (\ref{eq:remap_bnd}), (\ref{eq:remap_mass})  and (\ref{eq:remap_energy}) form a system of equations with $N$+4 equations and $N$ unknowns. We write the system of equations in matrix form as follows;
\begin{align} \label{matrix_a}
\begin{bmatrix}
      1 	&0 		&\dots 		&0 \\
      0 	&1 		&\dots 		&0 \\
      \vdots         			   \\
      0 	&0 		&\dots 		&1 \\
	h_{1/2}(\xi_{b=2}) &h_{3/2}(\xi_{b=2}) &\dots     &h_{N-1/2}(\xi_{b=2}) \\
	h_{1/2}(\xi_{b=1}) &h_{3/2}(\xi_{b=1}) &\dots     &h_{N-1/2}(\xi_{b=1}) \\
	w_{1/2}  &w_{3/2} &\dots &w_{N-1/2}   \\
	A_{1/2}  & A_{3/2} &\dots & A_{N-1/2} 
\end{bmatrix}
\begin{bmatrix}
	{\phi}^{n+1}(\xi^p_{1/2}) \\
	{\phi}^{n+1}(\xi^p_{3/2}) \\
\vdots \\
{\phi}^{n+1}(\xi^p_{N-1/2}) 
\end{bmatrix} =
\begin{bmatrix}
	\widehat{\phi}^{n+1}(\xi_{1/2}^p) \\
	\widehat{\phi}^{n+1}(\xi_{3/2}^p) \\
\vdots         \\
	\widehat{\phi}^{n+1}(\xi_{N-1/2}^p)  \\
	\widehat{\phi}^{n+1}_{b\star=1}  \\
	\widehat{\phi}^{n+1}_{b\star=2}   \\
	M^{n} + \Delta t  \mathcal{F}_m^n \\
	E^{n} + \Delta t  \mathcal{F}_e^n \\
\end{bmatrix}
\end{align}
Here, we use the notation $A_{i+1/2}=w_{i+1/2}\phi^{k-1}(\xi^p_{i+1/2})$ to prevent clutter in (\ref{matrix_a}).
If the energy constraint is enforced, then an
iterative solution to the system is necessary and
the superscript on the $\hat{\phi}$ should be $k$ rather
than $n+1$.  The counter $k$ is the iteration number.
If the energy constraint is not enforced, then the last row
disappears and the system can be solved directly.
A Least squares fit is used to solve the over determined system to obtain ${\phi}^{n+1}(\xi_j^p)$.

\subsubsection*{Two dimensional scheme}
\begin{figure}
	\centering 
	\mbox{
		\includegraphics[width=0.9\textwidth]{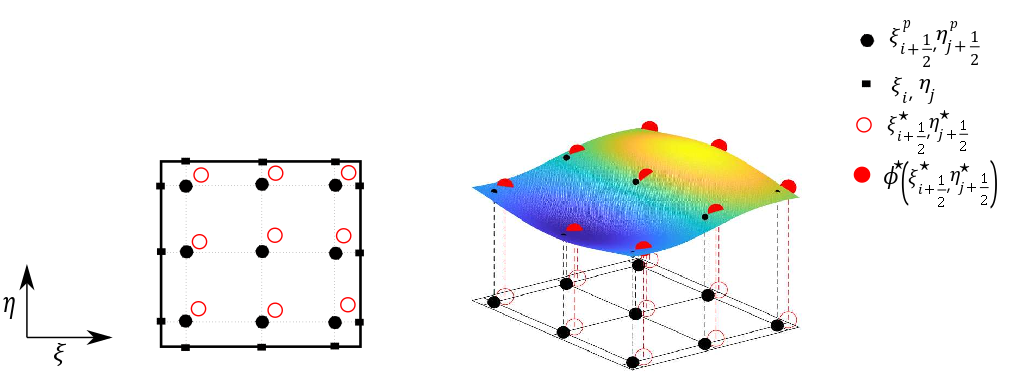}
	}
	\caption{Schematic of the DSEM SL method in two dimensions. The solution,$\phi^n$ is initialized on the Gauss nodes, $\xi^p_{i+1/2},\eta^p_{j+1/2}$. A 2D interpolant constructed at the Gauss nodes to obtain the intermediate solution using the advected solution, $\phi^\star$ on the advected positions, $\xi^\star_{i+1/2},\eta^\star_{j+1/2}$. }  
	\label{fig:2d_method}
\end{figure} 
Similar to the one dimensional scheme, in two dimensions the particles are initialized at the Gauss quadrature points, but along tensorial grid lines (Fig.\ \ref{fig:2d_method}), 
\begin{equation} \label{eq:5.1}
\phi^n(\xi,\eta) = \sum \limits_{i=0}^{N-1} \sum\limits_{j=0}^{N-1}{\phi}^n(\xi^p_{i+1/2},\eta^p_{j+1/2}) H_{i+1/2,j+1/2}(\xi,\eta), 
\end{equation} 
where $H_{i+1/2 j+1/2}(\xi,\eta) =h_{i+1/2}(\xi) h_{j+1/2}(\eta)$ and $h_{i+1/2}(\xi)$ and $h_{j+1/2}(\eta)$ are the Lagrange polynomials of degree $N-1$ along the $\xi$ and $\eta$ directions respectively.
Forward time integration of the particles is performed in $x$- and $y$- directions using the velocities $u(\xi^p_{i+1/2},\eta^p_{j+1/2})$ and $v(\xi^p_{i+1/2},\eta^p_{j+1/2})$, respectively, as prescribed or interpolated from an Eulerian formulation.
The advected particle positions and solutions so obtained are denoted as  $\xi^\star_{i+1/2},\eta^\star_{j+1/2}$ and $\phi^\star(\xi^\star_{i+1/2},\eta^\star_{j+1/2})$, respectively.

Up to this point, the algorithm is exactly like the one-dimensional algorithm but along tensorial grid lines.
The interpolation in two dimensions is more intricate as compared to the one dimensional algorithm. Instead of using an interpolant based on advected nodal values, a two-dimensional interpolant defined at the tensor grid of Gauss quadrature points is used to determine the intermediate solution, $\widehat{\phi}^{n+1}$.
\begin{eqnarray} \label{eq:5.4}
\phi^\star(\xi^\star_{i+1/2},\eta^\star_{j+1/2}) = \sum \limits_{l=0}^{N-1} \sum\limits_{m=0}^{N-1}\widehat{\phi}^{n+1}(\xi^p_{l+1/2},\eta^p_{m+1/2})  H_{l+1/2, m+1/2}(\xi_{i+1/2}^\star,\eta_{j+1/2}^\star), \\ \qquad i,j =0,1,...,N-1. \nonumber
\end{eqnarray} 
The interpolation then comes down to determining $\widehat\phi^{n+1}$ 
by inverting the matrix defined by $H_{l+1/2,m+1/2}(\xi^\star_{i+1/2},\eta^\star_{j+1/2})$ at every time step, \\
\resizebox{0.92\textwidth}{!}{
$ \left[	
 \begin{array}{llcl}
	H_{1/2,1/2}(\xi_{1/2}^\star, \eta_{1/2}^\star)    &H_{1/2, 3/2}(\xi_{1/2}^\star, \eta_{1/2}^\star) &\dots     &H_{N+1/2, N+1/2}(\xi_{1/2}^\star, \eta_{1/2}^\star) \\
	H_{1/2,1/2}(\xi_{1/2}^\star, \eta_{3/2}^\star)    &H_{1/2, 3/2}(\xi_{1/2}^\star, \eta_{3/2}^\star) &\dots     &H_{N+1/2, N+1/2}(\xi_{1/2}^\star, \eta_{3/2}^\star) \\
	 \vdots         			   \\
    H_{1/2,1/2}(\xi_{N+1/2}^\star, \eta_{N+1/2}^\star)    &H_{1/2, 3/2}(\xi_{N+1/2}^\star, \eta_{N+1/2}^\star) &\dots     &H_{N+1/2, N+1/2}(\xi_{N+1/2}^\star, \eta_{N+1/2}^\star)	
\end{array} \right] 
$
}
\\

For the matrix to not be ill-posed the particles should not cross each other or stray to far from the quadrature points, on which we know that the interpolant is well-defined. To keep the particles close to the quadrature points, we restrict the time step as follows,
\begin{equation} \label{eq:2d_delt_condition}
\Delta t \leq \text{min} \left[ \frac{\Delta \xi_{\text{min}}}{|u|_{max}}, \frac{\Delta \eta_{\text{min}}}{|v|_{max}} \right].
\end{equation}
With this restriction particles cannot leave the subdomain. In the interior of the subdomain the spacing between the quadrature nodes is quite a bit larger and the particles remain relatively close to the quadrature nodes. 
Here, $\Delta \xi_{min}$ and $\Delta \eta_{min}$ are the minimum grid
spacing between two particles and $|u|_{max}$ and $|v|_{max}$
are the maximum advection speed along $\xi$ and $\eta$ directions, respectively. Like for  the one-dimensional scheme, this time step limit is typically less restrictive
than the explicit CFL time step for the Eulerian DSEM solver. 
Once the intermediate solution $\hat\phi^{n+1}$ is obtained, we use 1D Lagrange polynomials to interpolate at the interface and boundary points.
\begin{eqnarray} \label{eq:2dx_bnd_interpolate}
\hat{\phi}^{n+1}_{j}|_b = \sum \limits_{i=0}^{N-1}{\hat{\phi}^{n+1}}(\xi_{i+1/2},\eta_{j+1/2}) h_{i+1/2}(\xi_b)  \qquad j =0,1,..N-1 \\ \qquad b=1,2 \nonumber
\end{eqnarray} 
\begin{eqnarray} \label{eq:2dy_bnd_interpolate}
\hat{\phi}^{n+1}_{i}|_b= \sum \limits_{j=0}^{N-1}{\hat{\phi}^{n+1}}(\xi_{i+1/2},\eta_{j+1/2}) h_{j+1/2}(\eta_b)  \qquad i =0,1,..N-1 \\ \qquad b=3,4 \nonumber
\end{eqnarray} 
where $b=1,2$ are the interface and boundary points along $\xi$ direction and $b=3,4$ are the interface and boundary points along $\eta$ direction. The final step is to use a  least squares remapping.

\section{Numerical tests}
\label{sec:4}
We test the semi-Lagrangian scheme for  combinations of accuracy of the time integrator and mass and/or energy constraints as listed in Table \ref{tab:cases}. We consider a one-dimensional constant velocity advection, a non-constant velocity advection and a two-dimensional constant velocity advection case. The tests and their results are discussed in the sections below.
\begin{table}[h] 
	\centering
	\begin{tabular}{|l|c|c|c|} \hline 
		Case  & Time integration order & Mass constraint & Energy constraint \\ \hline
		Basecase1 &  First & - & - \\ 
		MF1      &  First & Included & -\\ 
		MEF1     &  First & Included & Included\\ 
		Basecase2 & Second & - & - \\
		MF2       &  Second & Included & -\\ 
		MEF2      &  Second & Included & Included \\ 
		MF3       & Third  & Included & - \\
		\hline 
	\end{tabular} 
	\caption{Definition of acronyms used to describe combinations of accuracy of the time integration and types of constraints. }
	\label{tab:cases}
\end{table}

In each of the tests, we inspect the behavior of local errors and  global (error) norms. We consider the $L_2$ error norm which is calculated by summing up local $L_2$ error norms in each subdomain, $k$, as,
\begin{equation} \label{eq:5.2}
\|e\|_{L^2} = \sum \limits_{k=1}^{K} \sqrt{\int_{\Omega_k} (\phi-\phi_{\text{exact}})^2 d\xi}.
\end{equation}
The conservation properties of the method are inspected  with the following global mass and energy norms, 
\begin{equation}
\|M\|=\sum \limits_{k=1}^{K} \frac{\int_{\Omega_k} \phi d\xi}{\int_{\Omega_k}\phi_\text{exact} d\xi},
\end{equation}
and
\begin{equation}
\|E\| = \sum  \limits_{k=1}^{K} \frac{\int_{\Omega_k}\phi^2 d\xi}{\int_{\Omega_k}\phi^2_\text{exact} d\xi}.
\end{equation}
respectively.
For the test cases in  which the analytic global mass is zero, we adjust the definition of the global mass norm as follows 
$\|M_r\|=\sum \limits_{k=1}^{K} \int_{\Omega_k} \phi d\xi$ to avoid division by zero.

We expect that the accuracy of the scheme depends on the interpolation and time integration accuracy and the impact of the constraints through the least squares fitting. We will discuss error and stability behavior by means of the test below.
\subsection{One dimensional constant velocity}

As a first test, we consider a linearly advected sine wave according to (\ref{eq:2.5}) and (\ref{eq:2.4}) with $u$=1   
in a domain $x$=$[0,1]$.  The initial condition is $\phi(x,0)$ =$\sin(2\pi x)$. Periodic boundary conditions are specified so that the sine wave propagates through the domain.  Simulations are run until $T=$10, at which time error norms are determined and compared for several variations of the numerical scheme according to Table\ \ref{tab:cases}. Tests are conducted for polynomial orders of $P$=4,5,6, and 7 with different number of elements, $H$=4,5,6, and 7.

\begin{figure} 
\centering 
\mbox{ 
  \includegraphics[width=0.5\textwidth]{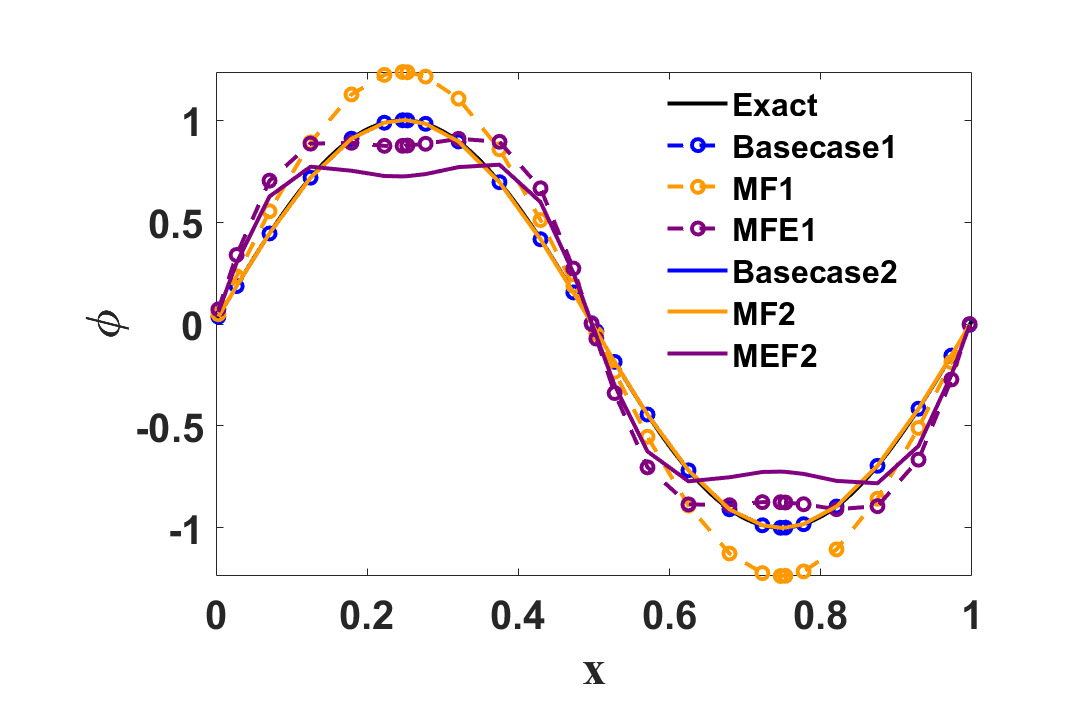} 
  \hspace{0.01\textwidth}
  \includegraphics[width=0.5\textwidth]{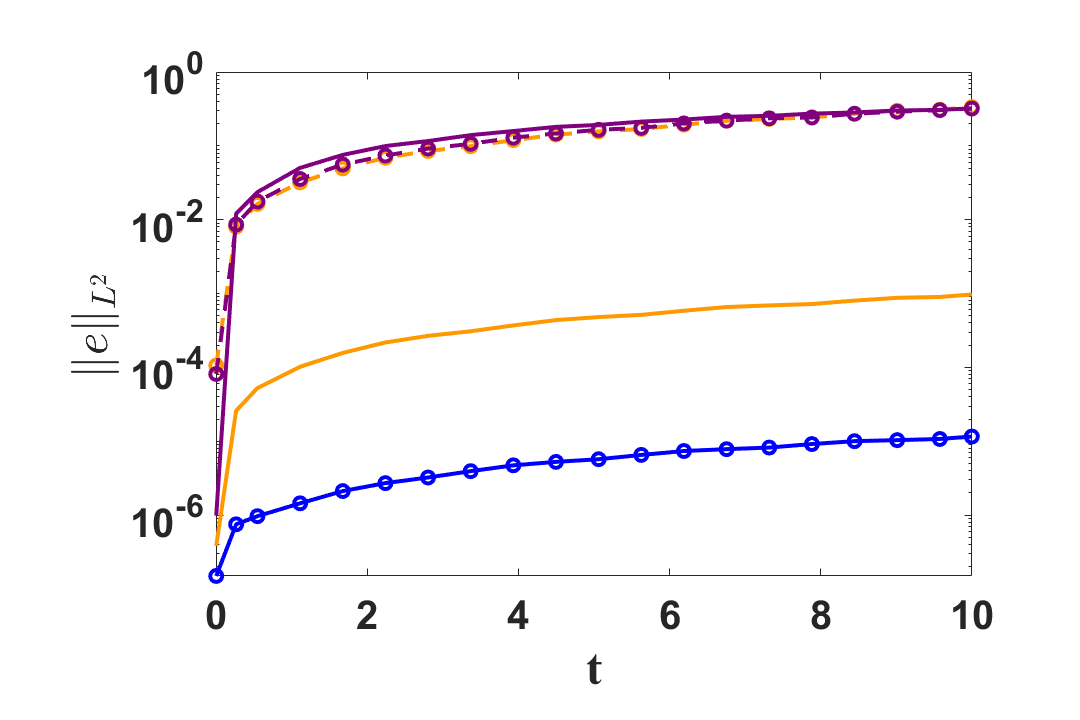}
}
\mbox{
\hspace{0.4cm}
\makebox[0.45\textwidth]{(a)}
 \hspace{0.05\textwidth}
\makebox[0.45\textwidth]{(b)}
}

\mbox{
  \vspace{1.5cm}
  \includegraphics[width=0.5\textwidth]{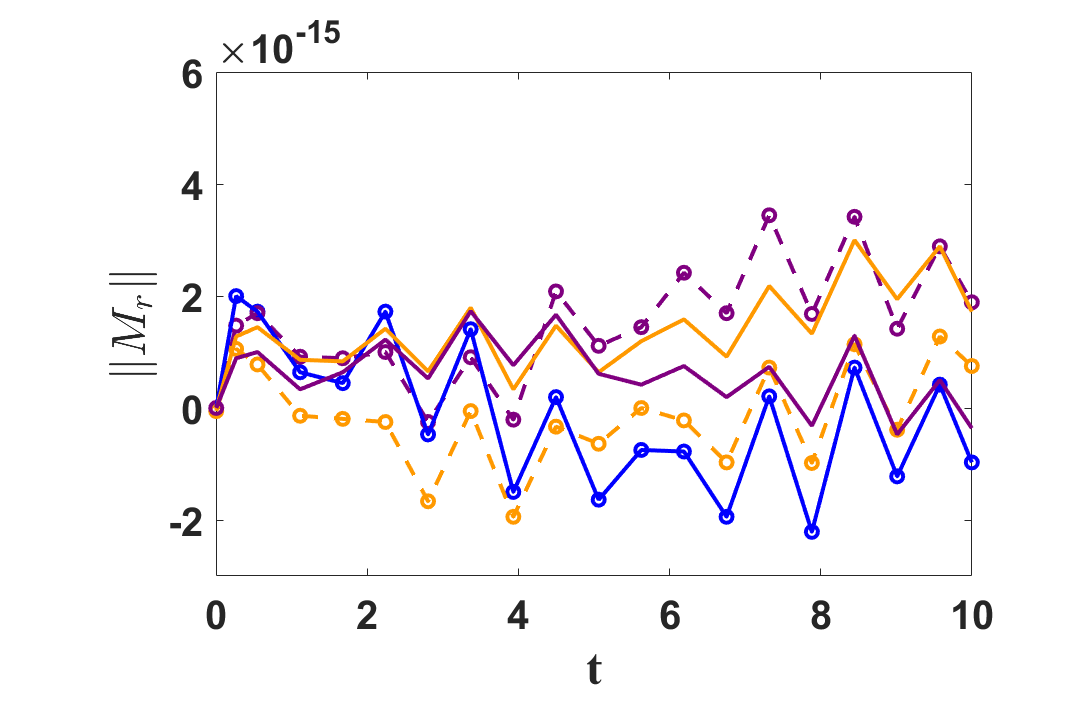} 
   \hspace{0.01\textwidth}
  \includegraphics[width=0.5\textwidth]{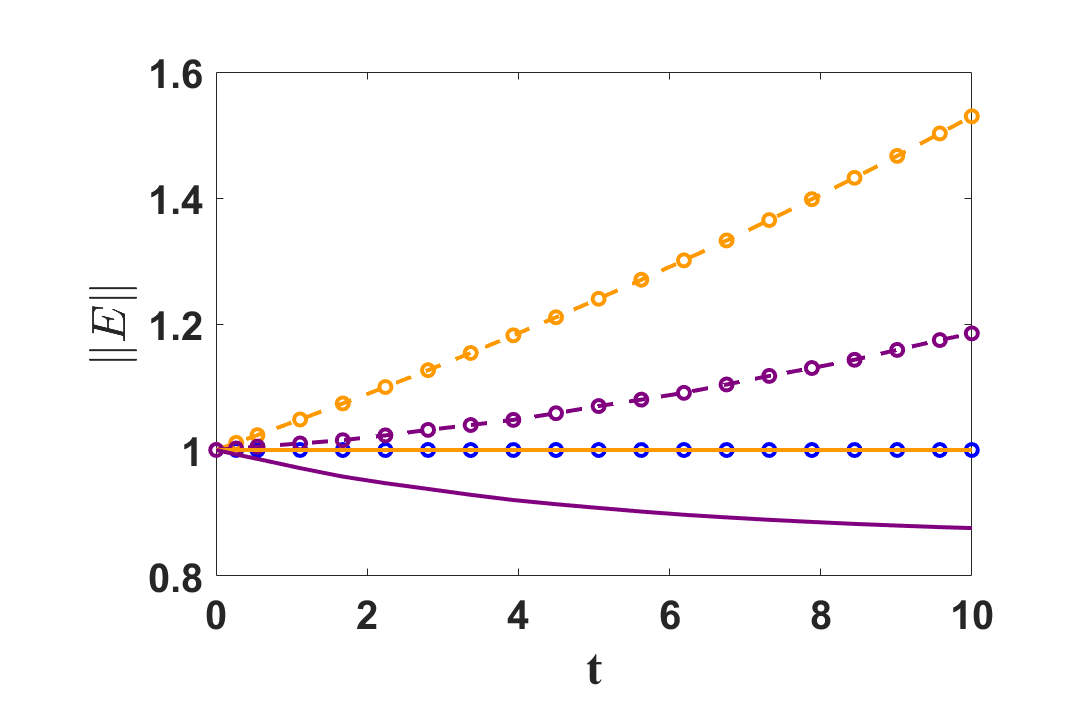}
} 
\mbox{
\hspace{0.4cm}
\makebox[0.45\textwidth]{(c)}
 \hspace{0.05\textwidth}
\makebox[0.45\textwidth]{(d)}
}
\caption{Results for a linearly advected sine wave with DSEM-SL using $H$=4 subdomains and an approximation order of $P$=6; (a) the solution, $\phi$, versus $x$  at $T$=10. 
The $L^2$ error norm, $\|e\|_{L^2}$, the mass norm, $\|M\|$, and the energy norm, $\|E\|$ are plotted versus time, $t$, in subfigures (b), (c) and (d), respectively.}  
\label{fig:1_1}
\end{figure}

 
Because the advection velocity is constant, the time integration of the particle's location is exact, independent of the accuracy of the time integration. The accuracy of the Basecases (Table\ \ref{tab:cases}) which only use the boundary constraint in the least squares fit step of the algorithm do not depend on the time integrator. So, the time evolution of the global norms for the Basecases as shown in Figure\ \ref{fig:1_1}b, c and d for a case with $H$=4 and $P$=6, do not depend on the accuracy of the time integrator. 
We hence also observe that the results at $T$=10 of Basecase1 and Basecase2 in Figure \ref{fig:1_1}a overlap each other.

The temporal update of mass and energy
in (\ref{eq:remap_mass}) and (\ref{eq:remap_energy}) used as  constraints in (\ref{matrix_a})  {\em do} depend on the accuracy of the time integrator.
The first order schemes, which can be expected to introduce a time integration error  on the order of the explicit stable time step, which is  $\Delta t ={O}(10^{-3})$ for this case, yield a global $L^2$ error at $T$=10 of $\mathcal{O}(10^{-1})$ for both mass and energy constraint cases (Fig.\ \ref{fig:1_1}b) .
While the mass norms are constant (Fig.\ \ref{fig:1_1}c), 
the  energy norms (Fig.\ \ref{fig:1_1} d) are not conserved and, in fact, increase over time. 
This is further confirmed by the increased amplitude of the sine wave in the solution at $T$=10 in Figure\ \ref{fig:1_1}a. Clearly, the local time error accumulates in time and the first order temporal update of the constraint in combination with the least squares fit is leading to a significant increased error as compared to the Basecases.Using the second order time integration, improves the $L^2$ error for the MF2 cases, i.e. the mass flux constraint.  Moreover, for this case the global energy conservation is accurate. 

\begin{figure} 
\centering 
\mbox{ 
  \includegraphics[width=0.5\textwidth]{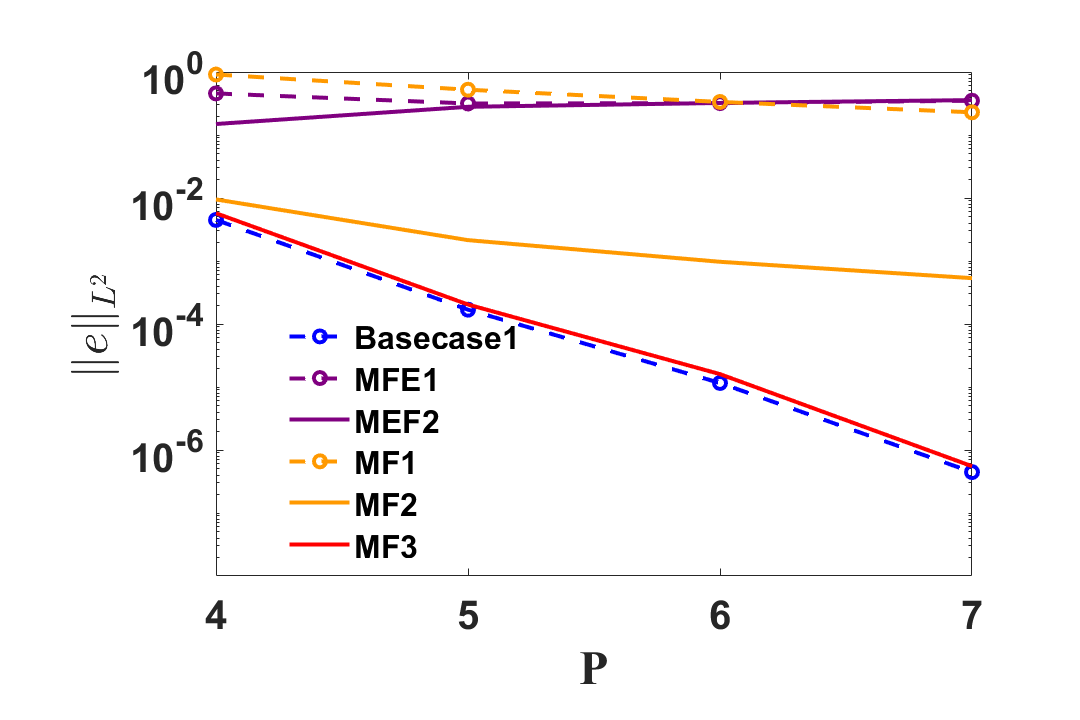} 
}
\caption{$P$-convergence of the  schemes in Table\ \ref{tab:cases} for a constant velocity linearly advected sine wave using $H$=4 subdomains  at $T$=10 using the explicit stable time step.} 
\label{fig:1_2}
\end{figure}

Enforcing the additional constraint on the local energy (MEF2) produces spurious modes in the solution as observed in cases MEF1 and MEF2 shown in Figure \ref{fig:1_1}a. It would appear that the extra energy constraint leads to the least-squares system being overconstrained, in turn, leading to instability. This is underscored by the non-conservation of the global energy for the MEF cases. 

The $\|e\|_{L^2}$ plotted versus $P$ for $T$=10 in Figure \ref{fig:1_2} shows that the Base cases are exponentially convergent. While the mass and energy is not formally conserved for those cases, the high order convergence and accuracy  translates to an accurate, spectral approximation of the mass and energy. The first-order MF1 case does not show $P$-convergence in Figure\ \ref{fig:1_2} because the overall error is bound by the accuracy of the time integrator. The second-order MF2, which is more accurate in time, does show an error improvement with $P$ that is algebraic. This is consistent with the algebraic convergence of the time integration and the stable time step reduction that is inherent to an increase of $P$. 
With a further increase of the time accuracy to third-order, the MF3 case matches the result of the Basecase, which is bound by the accuracy of the polynomial interpolation.
Table, \ref{tab:h_conv} shows that the Basecase methods exhibit formal  $H$-convergence. The order of convergence is a little lower than $P$ which is the expected convergence rate. The MF1 and MF2 cases converge according to the order of time integration. 
%

\begin{table}
\centering
\begin{tabular}{|c|c|c|c|} \hline
\multicolumn{2}{|c|}{} & \multicolumn{2}{|c|}{$P$=4} \\ \cline{3-4}
\multicolumn{2}{|c|}{} & $\|e\|_{L^2}$ & Slope   \\ \hline
Basecase & $H$=5 & $2.067 \times 10^{-3}$ & 3.47 \\
         & $H$=6 & $1.07  \times 10^{-3}$ & 3.60 \\
         & $H$=7 & $5.98  \times 10^{-4}$ & 3.78  \\ \hline
 MF1     & $H$=5 & $8.02  \times 10^{-1}$ & 0.57  \\
         & $H$=6 & $7.38  \times 10^{-1}$ & 0.45  \\
         & $H$=7 & $6.65 \times 10^{-1}$  & 0.67  \\ \hline
 MF2     & $H$=5 & $5.53 \times 10^{-3}$   & 2.41 \\
         & $H$=6 & $3.69 \times 10^{-3}$   & 2.22  \\
         & $H$=7 & $2.737 \times 10^{-3}$  & 1.93  \\ \hline
MF3     & $H$=5  & $2.64  \times 10^{-3}$  & 3.47 \\
        & $H$=6  & $1.39  \times 10^{-3}$  & 3.52  \\
        & $H$=7 &  $7.92  \times 10^{-4}$  & 3.65   \\ \hline         
\end{tabular} 
 \caption{$H$ convergence of the $L^2$ error norm, $\|e\|_{L^2}$, for $P$=4 at $T$=10.}
\label{tab:h_conv}
\end{table}

\begin{figure}
\centering 
\begin{subfigure}{.47\textwidth}
  \includegraphics[width=1\linewidth]{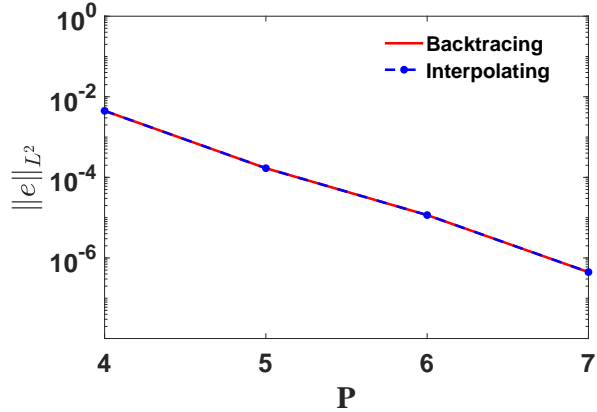}
\end{subfigure}%
\caption{Simulation of constant velocity advection using DSEM-SL method using the explicit stable time step. Plot comparing the backtracking and interpolation boundary constraints when $H=4$ and $t=10$.} 
\label{fig:bc}
\end{figure}

In Figure \ref{fig:bc} we compare solutions obtained with two approaches to determine the interface solution that include  the backtracking approach and the interpolation approach. These interface values  are applied as the boundary constraint in the least squares fit. The $L^2$ error with $H$=4 with $P$=4,5,6, and7 shows that both methods show no distinguishable result. We find that the interpolation error is typically small and that it matches the solution from backtracking. The interpolation approach is easier to implement and we therefore prefer it.
 
\begin{figure} 
\centering 
\mbox{ 
  \includegraphics[width=0.5\textwidth]{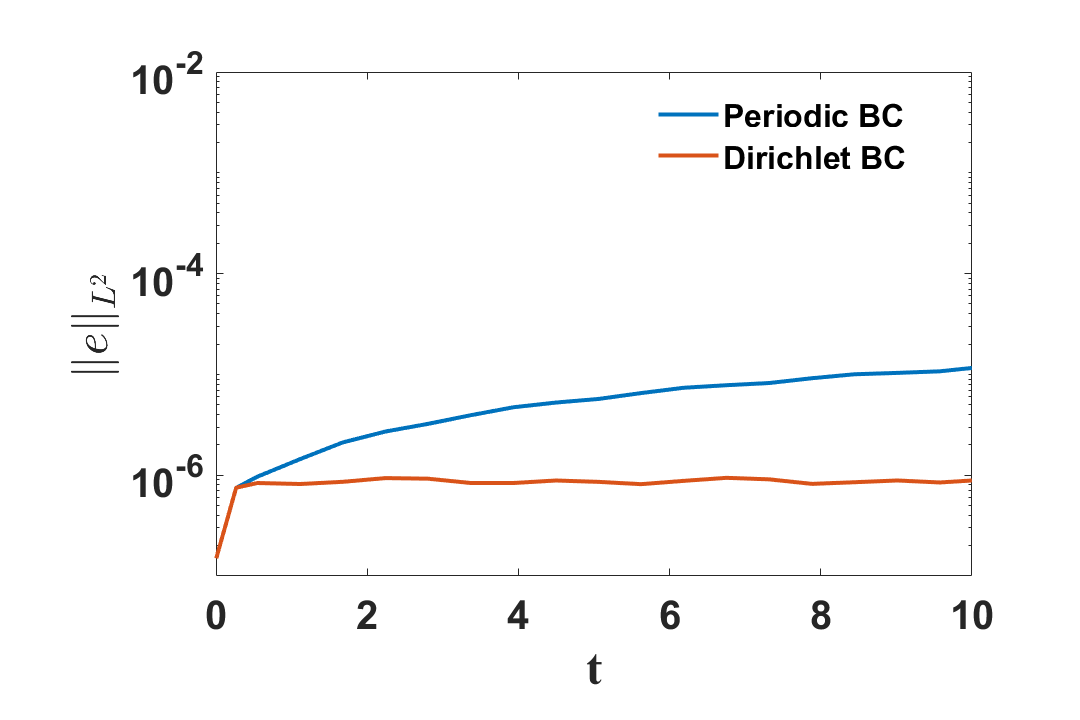} 
   \hspace{0.01\textwidth}
   \includegraphics[width=0.5\textwidth]{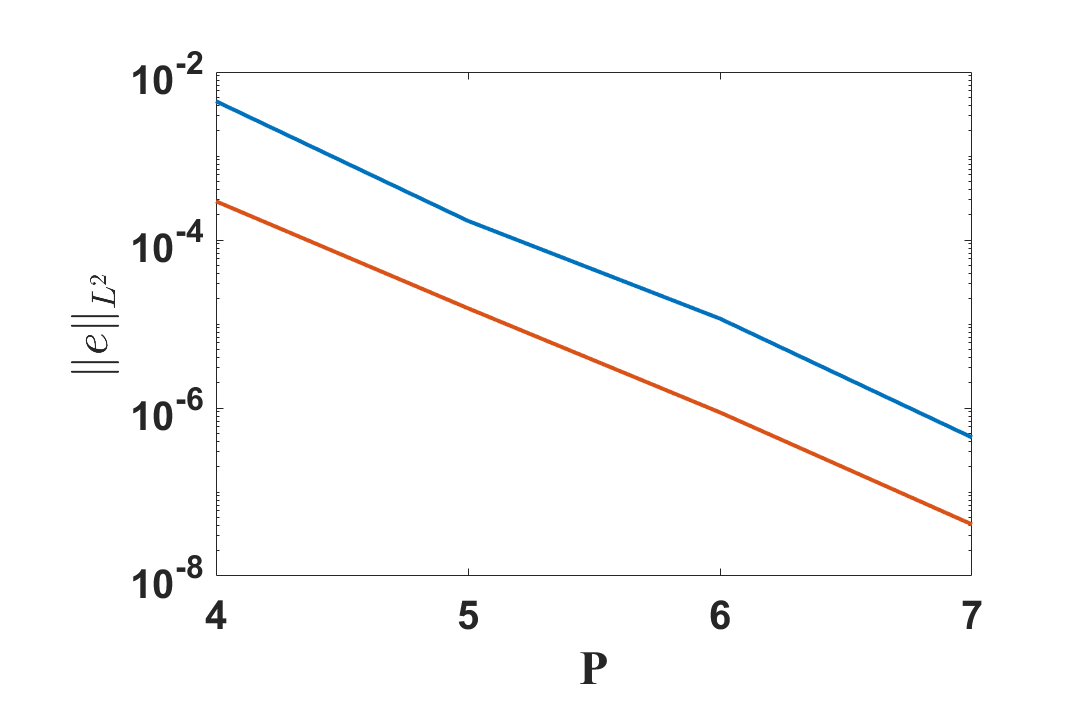}
}
\mbox{
 \hspace{0.4cm}
 \makebox[0.45\textwidth]{(a)}
  \hspace{0.05\textwidth}
 \makebox[0.45\textwidth]{(b)}
 }
\caption{Comparing Dirichlet and periodic boundary conditions using the Basecase for a constant velocity linearly advected sine wave. (a) shows the time evolution of the $\|e\|_{L^{2}}$ error using $P$=6 and $H$=4. (b) compares the $P$-comnvergence at $T$=10d.} 
\label{fig:1_3}
\end{figure}
As a last verification of the algorithm for the constant velocity case, we compare results obtained with a Dirichlet boundary condition those obtained with a periodic boundary condition that we have used thus far.  Figure \ \ref{fig:1_3}, shows that the Dirichlet boundary condition is slightly more accurate than the periodic boundary condition for the Basecase. This can be expected because the Dirichlet boundary condition is exact at the inflow, and from there errors accumulate/increase in space, whereas the periodic conditions do not have this exact reference and hence the errors are homogeneous in space.

\subsection{One dimensional non-constant velocity}

To study the effect of a non-constant advection velocity,
we consider a problem that is commonly considered
for semi-Lagrangian methods \cite{dgsl4} and  that solves
the  following conservation equation in Eulerian form
\begin{equation}
\frac{\partial \phi}{\partial t} + \frac{\partial (-\sin(x) \phi)}{\partial x} =0
\label{eq:usin_euler_conservative}
\end{equation}
on the domain $x \in [0, 2 \pi]$. The initial condition is
\begin{equation}
  \phi(x,0)= { {\sin \left(2 \tan^{-1} \left[\exp(-1) \tan({x \over 2})\right] \right)} \over {\sin(x)} }.  
\end{equation} 
Periodic boundary conditions are used. The problem has an analytical solution which is described by a symmetric distribution function which dampens into a constant valued solution over time.

We can rewrite (\ref{eq:usin_euler_conservative}) in non-conservation form as
\begin{equation}
\frac{\partial \phi}{\partial t} - \sin(x)\frac{\partial \phi}{\partial x} =\cos(x) \phi.
\label{eq:usin_euler_nonconservative}
\end{equation}
The equivalent Lagrangian form that we use as 
the governing equation for solution with the DSEM-SL method (following (\ref{eq:2.5}) and (\ref{eq:2.4})), is
\begin{eqnarray} \label{eq:usin_lagr}
\frac{d \mathbf{x^p}(t)}{d t} &=& -\sin(x) \nonumber \\
\frac{D \phi}{D t} &=& -\cos(x)\phi 
\end{eqnarray}  
Hence, the advection velocity,  $u$=$\sin(x)$, is non-constant.
We note that the right-hand side of the $\phi$ equation is
non-zero and will require a temporal update according to
(\ref{eq:4.8}). For the constant velocity case discussed in the previous
this section right hand side was zero, and the $\phi$ solution for that
case is simply constant along its characteristic path.

\begin{figure} 
\centering 
\mbox{ 
  \includegraphics[width=0.5\textwidth]{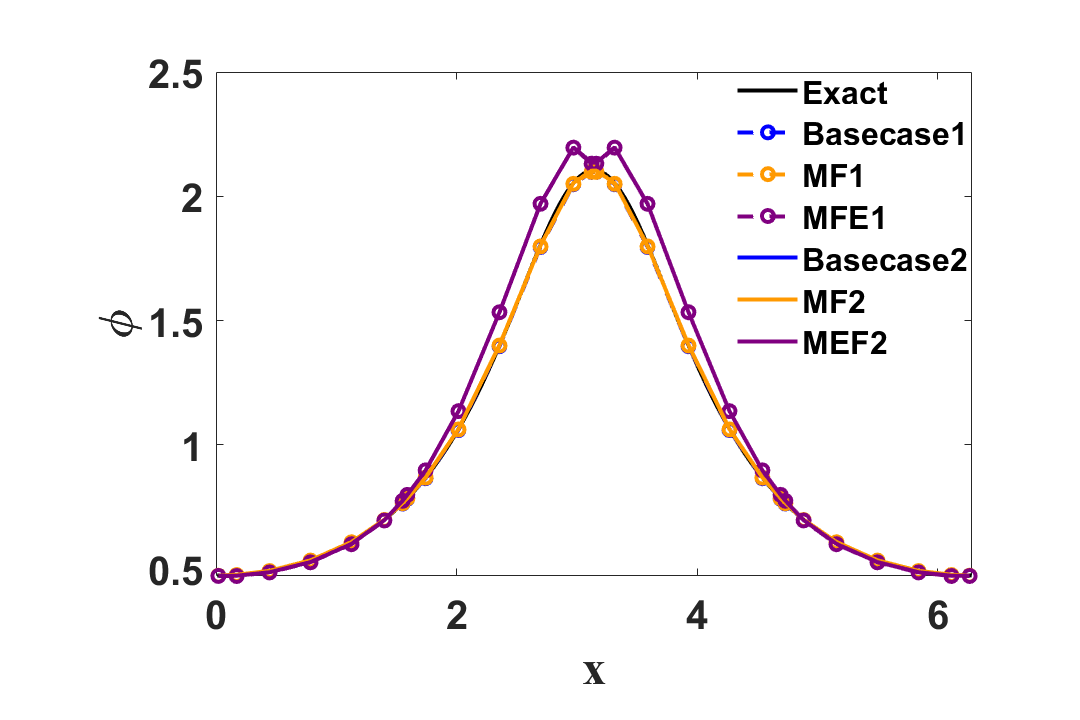} 
  \hspace{0.01\textwidth}
  \includegraphics[width=0.5\textwidth]{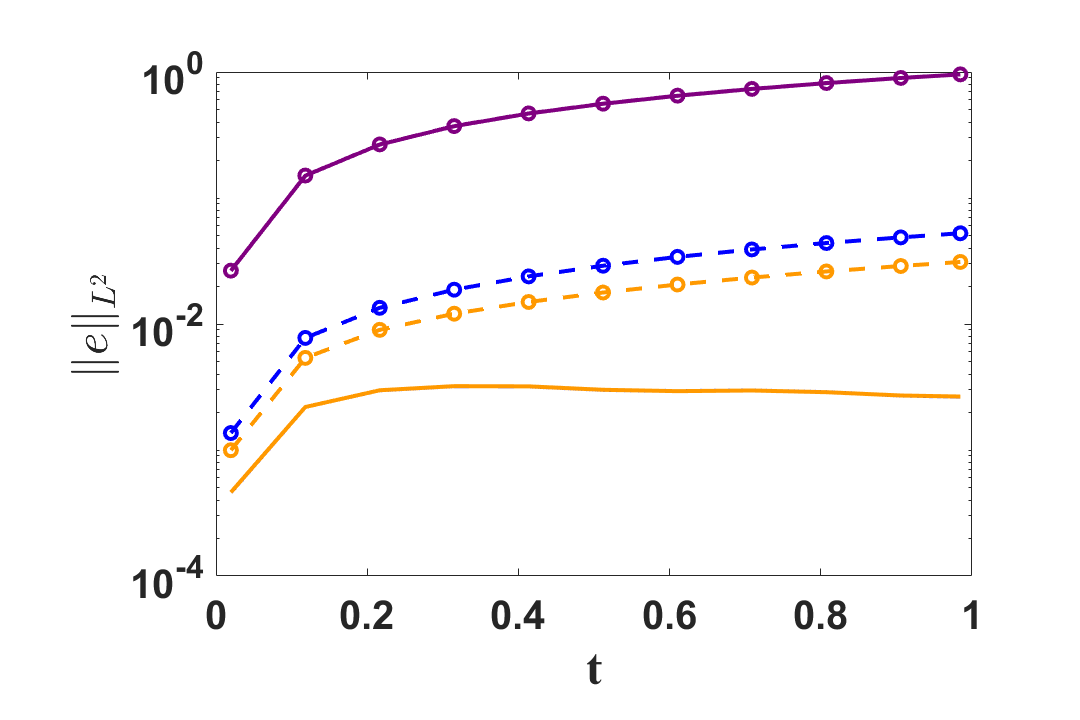}
}
\mbox{
\hspace{0.4cm}
\makebox[0.45\textwidth]{(a)}
 \hspace{0.05\textwidth}
\makebox[0.45\textwidth]{(b)}
}

\mbox{
  \vspace{1.5cm}
  \includegraphics[width=0.5\textwidth]{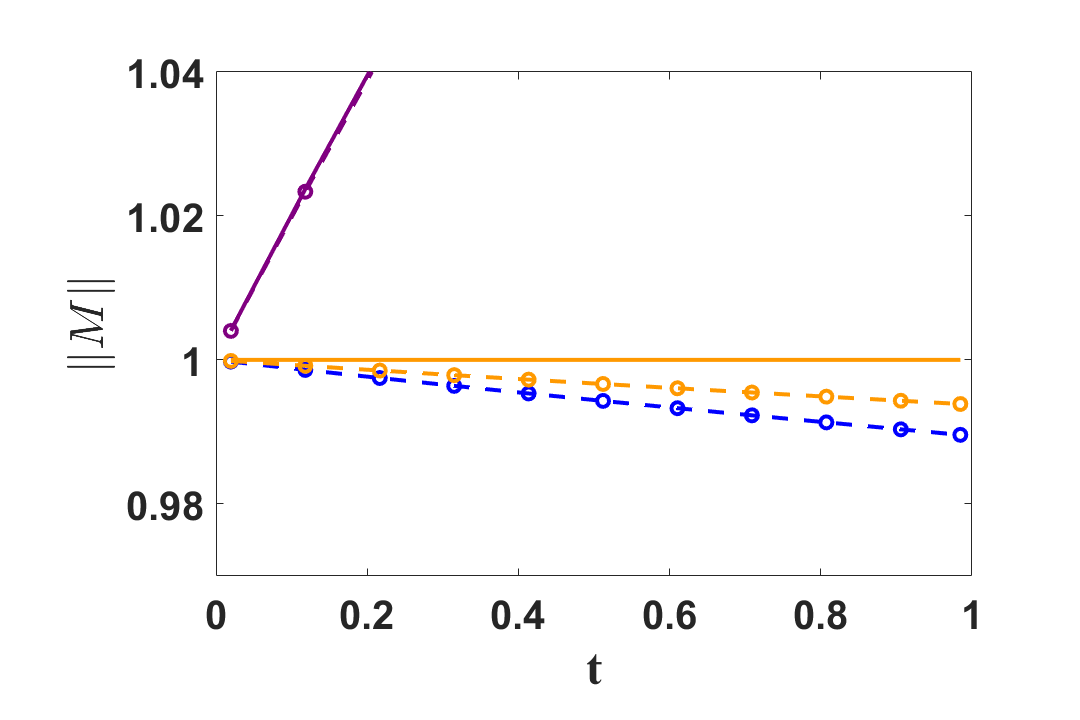} 
   \hspace{0.01\textwidth}
  \includegraphics[width=0.5\textwidth]{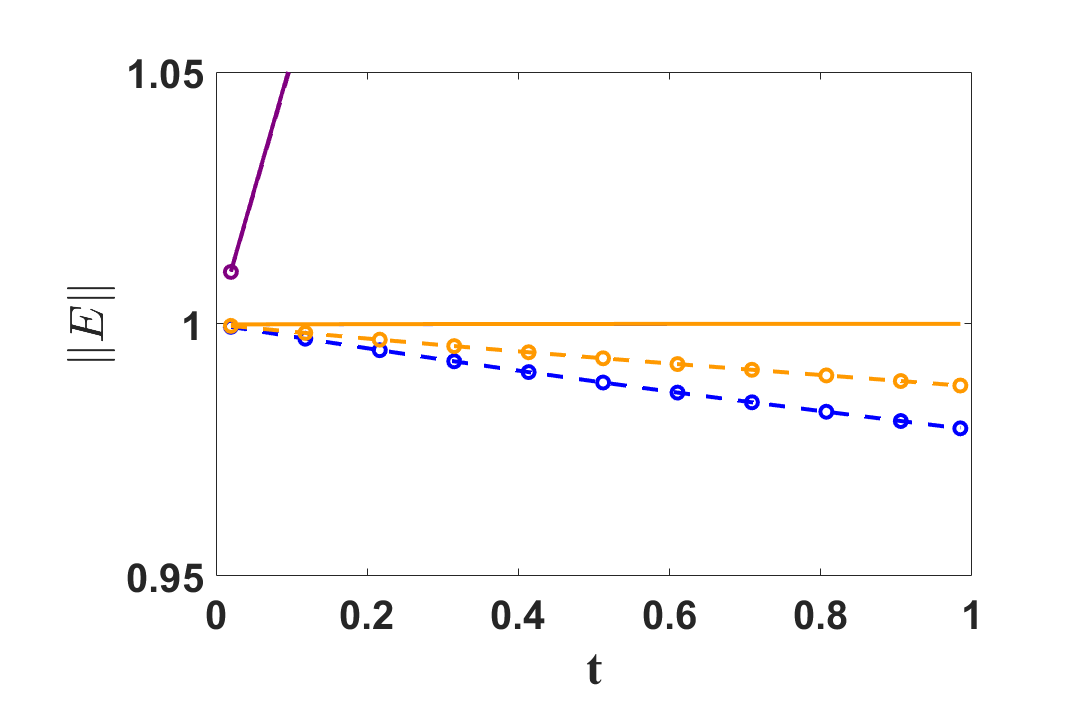}
} 
\mbox{
\hspace{0.4cm}
\makebox[0.45\textwidth]{(c)}
 \hspace{0.05\textwidth}
\makebox[0.45\textwidth]{(d)}
}
\caption{Results for a non-constant velocity advection with DSEM-SL using $H$=4 subdomains and an approximation order of $P$=6; (a) the solution, $\phi$, versus $x$  at $T$=1. 
The $L^2$ error norm, $\|e\|_{L^2}$, the mass norm, $\|M\|$, and the energy norm, $\|E\|$ are plotted versus time, $t$, in subfigures (b), (c) and (d), respectively.}  
\label{fig:2_1}
\end{figure}
Figure \ref{fig:2_1} shows the results for a case with four subdomains and $P$=6. The explicit time step
is $\Delta t$ is $1.97 \times 10^{-2}$. As opposed to the constant velocity case in the previous section, the $\phi$ solutions of Basecase with non-constant velocity are dependent on the time integrator because of the non-zero right hand side in the $\phi$ equation in (\ref{eq:usin_lagr}). The Basecase1 with a first order time integration has an $L^2$ error on the order of the time step of $\mathcal{O}(10^{-2})$ and is clearly bounded by the time integration error. The second order Basecase2 improves the accuracy to $\mathcal{O}(10^{-3})$ consistent with the increased time accuracy order. Because of the reduced accuracy for the first-order Basecase1, the global mass norm  reduces by $2\%$ over the considered time interval of $T$=1 (See also Table\ \ref{tab:3}). The second order Basecase2 conserves mass up to four decimals over $T$=1. 

\begin{table}[h] 
\centering
\begin{tabular}{|c|c|c|c|} \hline 
  Case  & $\|M\|$-1 & $\|E\|$ & $ \|e \|_{L^2}$    \\ \hline
 Basecase1 &  -0.0104  & 0.9793 & $5.25\times 10^{-2}$\\
 MF1      &   -0.0061  & 0.9878 & $3.11\times 10^{-2}$\\
 MEF1     &    0.1334  & 1.3148 & $9.57\times 10^{-1}$\\
 Basecase2&    0.0000  & 1.0000 & $2.65\times 10^{-3}$\\
 MF2      &    0.0000  & 1.0000 & $2.65\times 10^{-3}$\\
 MEF2     &    0.1367  & 1.3217 & $9.56\times 10^{-1}$\\
 MF3      &    0.0000  & 1.0000 & $2.63\times 10^{-3}$\\
 \hline
	\end{tabular} 
	\caption{Comparison of global for the non-constant advection velocity case at $T$=1 with $H$=4 and $P$=6.}
	\label{tab:3} 
\end{table}

The addition of a local mass constraint to the first order scheme, MF1 improves the global mass and energy conservation properties of the Basecase11 by $1\%$. The second order MF2 which includes the local mass constraint matches the solution of Basecase2. This is because the second order schemes reduce the time integration error to be on the order of  $\mathcal{O}(10^{-4})$ and because the error is dominated by the interpolation error. 

The addition of local energy constraints  produces spurious modes which are observed in Figure\ \ref{fig:2_1}a. This is similar to the constant advection velocity case. The overconstrained system is unstable as confirmed by the increase in relative global energy norm. 
 
\begin{figure} 
\centering 
\mbox{ 
  \includegraphics[width=0.5\textwidth]{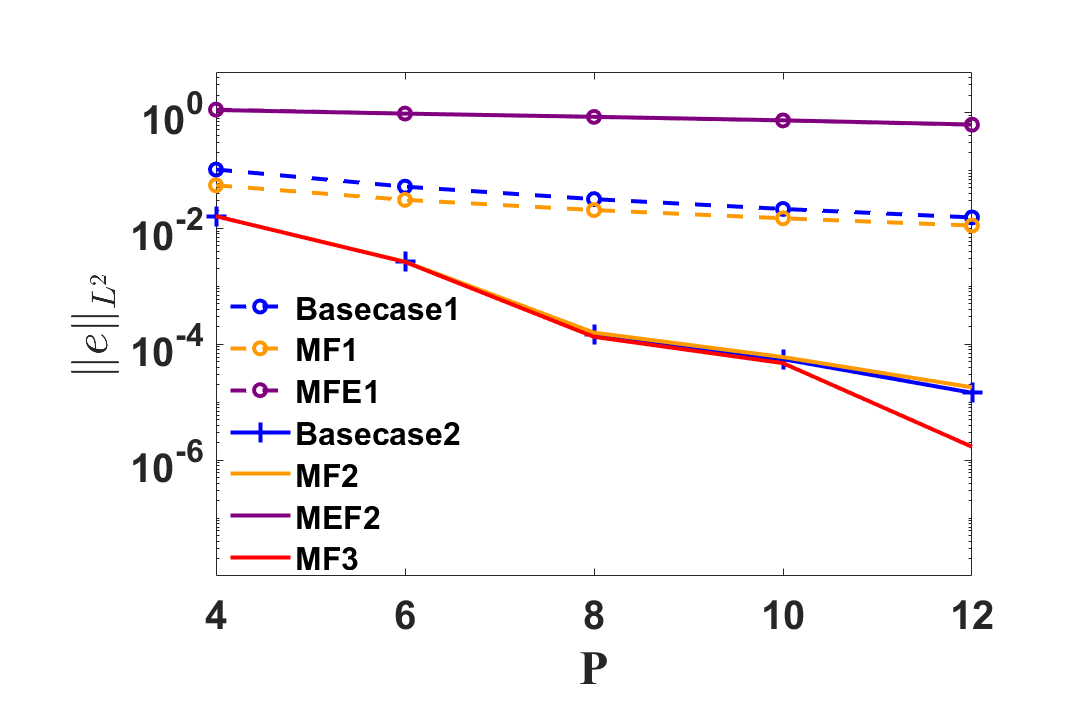} 
  \hspace{0.01\textwidth}
  \includegraphics[width=0.5\textwidth]{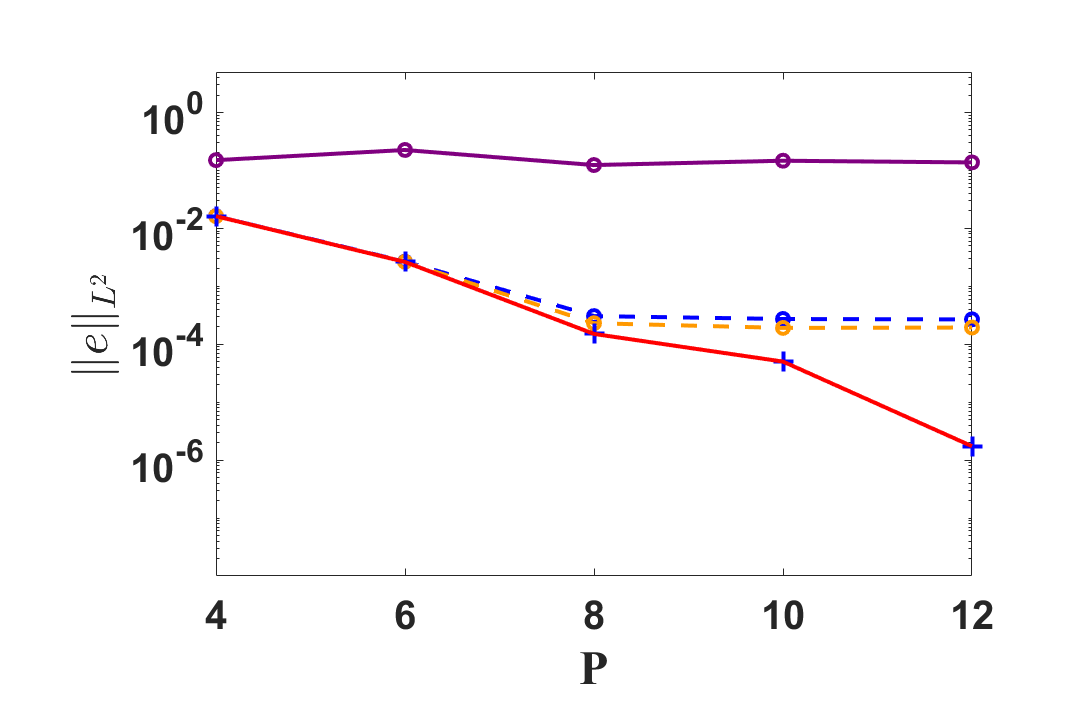}
}
\mbox{
\hspace{0.4cm}
\makebox[0.45\textwidth]{(a)}
 \hspace{0.05\textwidth}
\makebox[0.45\textwidth]{(b)}
}
\caption{$P$-convergence of the  schemes in Table\ \ref{tab:cases} for a non-constant velocity advection case using $H$=4 subdomains at $T$=1. for the explicit stable time step $\Delta t$ (a) and for a fixed time step $\Delta t$=$10^{-4}$(b) } 
\label{fig:2_2}
\end{figure}

Figure\ \ref{fig:2_2}a, which plots the approximation order $P$ against the $L^2$ error norm  shows mostly algebraic convergence for the first order Basecase1 and MF1 methods. The second order methods Basecase2 and MF2 show spectral convergence in $P$ until $P$=10, when the time integration error bounds the error. The third-order, MF3 converges exponentially till $P$=12.

To understand the effect of time step, we consider a smaller time step $\Delta t$ = $10^{-4}$. Figure \ref{fig:2_2}b  shows that this indeed improves the convergence of the first order methods, and spectral convergence can be  observed for this case up to $P$=8 where the error is of the order of $\mathcal{O}(10^{-4})$. Beyond $P$=8, the first order schemes are again bounded by the time integration error while the second and third order schemes continue to show spectral convergence. 

\subsection{Two dimensional constant velocity}
\noindent To test the algorithm in two dimensions,
we consider the two-dimensional constant velocity, linear advection equation in Lagrangian form
\begin{eqnarray}  \label{eq:5.3}
{d x \over d t} &=& a \nonumber \\
{d y \over d t} &=& b \nonumber \\
{d \phi \over dt} &=& 0
\end{eqnarray}
We solve this equation on the domain $x$=[0,1],$y$=[0,1] with $a$=2, and $b$=1. An initial condition is set as $\phi(x,y,0) = \sin(2\pi x)\sin(2\pi y)$. Periodic boundary conditions are applied in both $x$- and $y$-directions. The simulation is carried out until time $T=1$. 

Basecase1, Basecase2, MF1 and MF2 as listed in Table \ref{tab:cases} are simulated for a number of subdomains in $x$-
 and $y$-directions of $H_x$=$H_y$=4,5,6, and 7 and  approximation orders of $P$=4,5,6, and 7. 
 
\begin{figure} 
\centering 
\mbox{ 
  \includegraphics[width=0.5\textwidth]{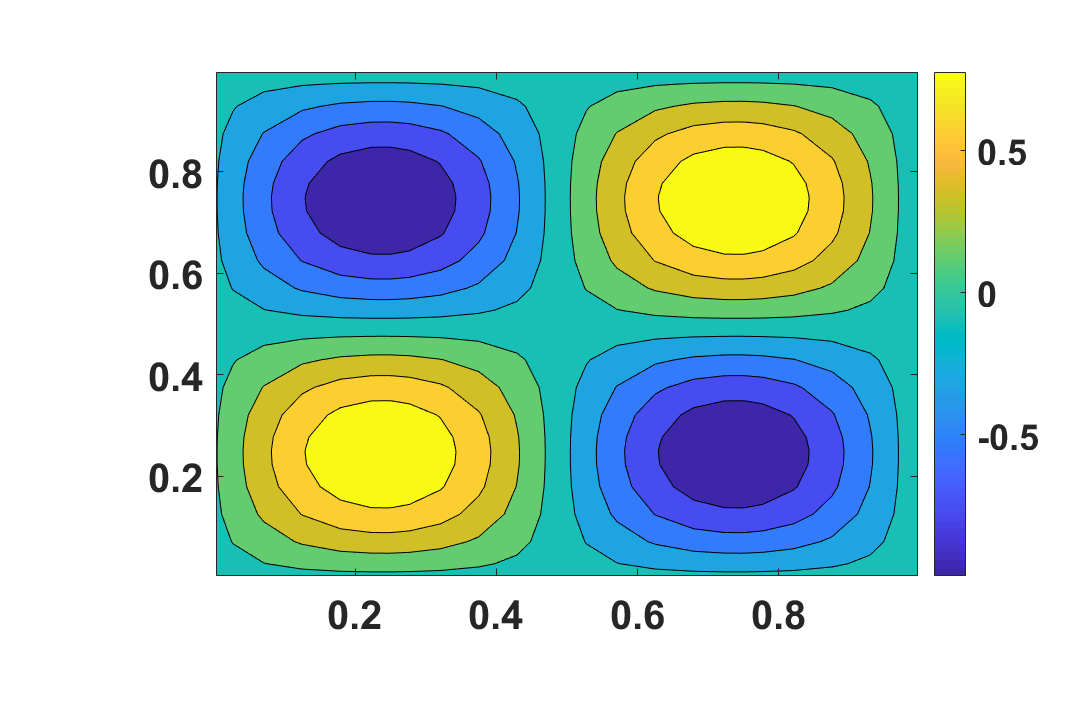} 
  \hspace{0.01\textwidth}
  \includegraphics[width=0.5\textwidth]{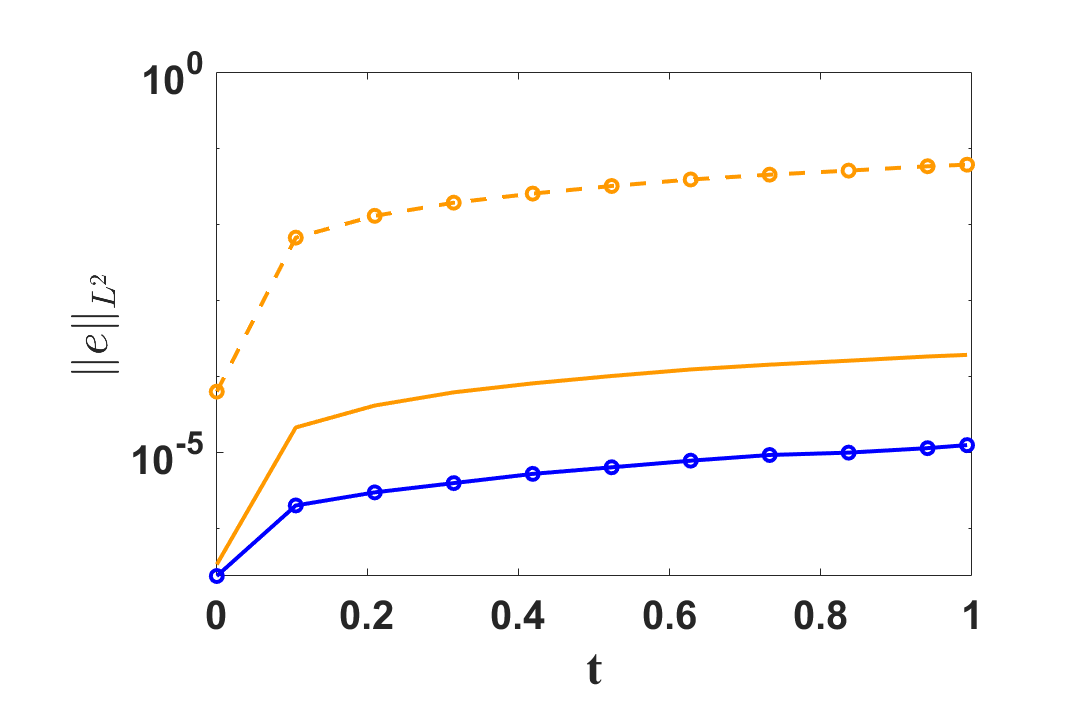}
}
\mbox{
\hspace{0.4cm}
\makebox[0.45\textwidth]{(a)}
 \hspace{0.05\textwidth}
\makebox[0.45\textwidth]{(b)}
}

\mbox{
  \vspace{1.5cm}
  \includegraphics[width=0.5\textwidth]{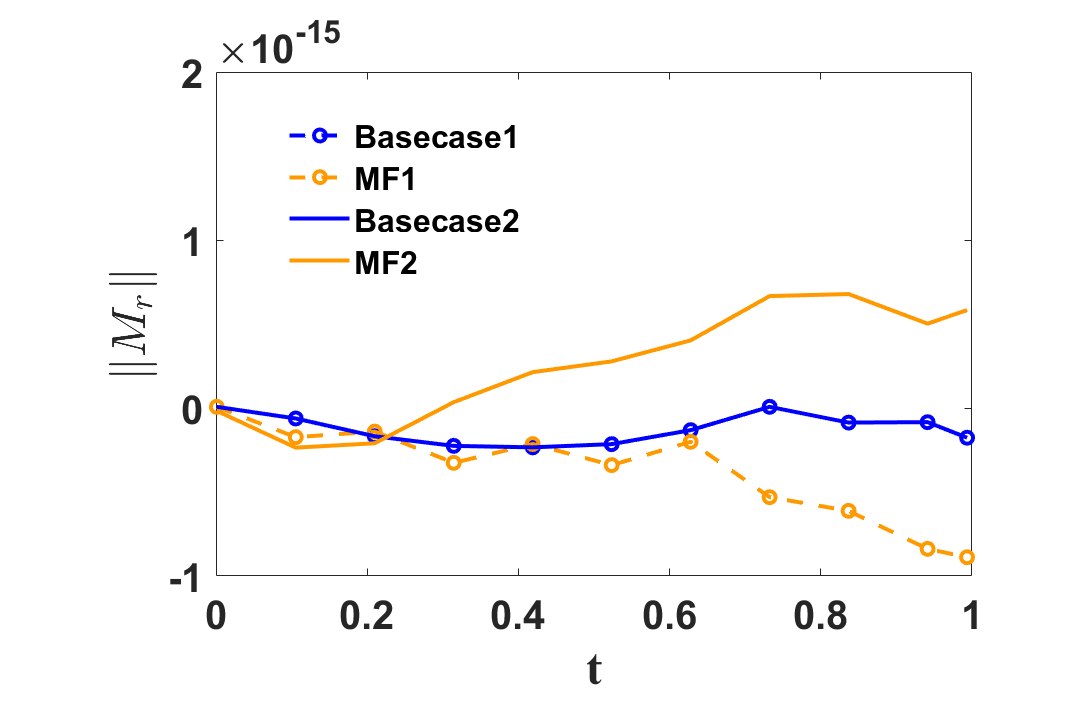} 
   \hspace{0.01\textwidth}
  \includegraphics[width=0.5\textwidth]{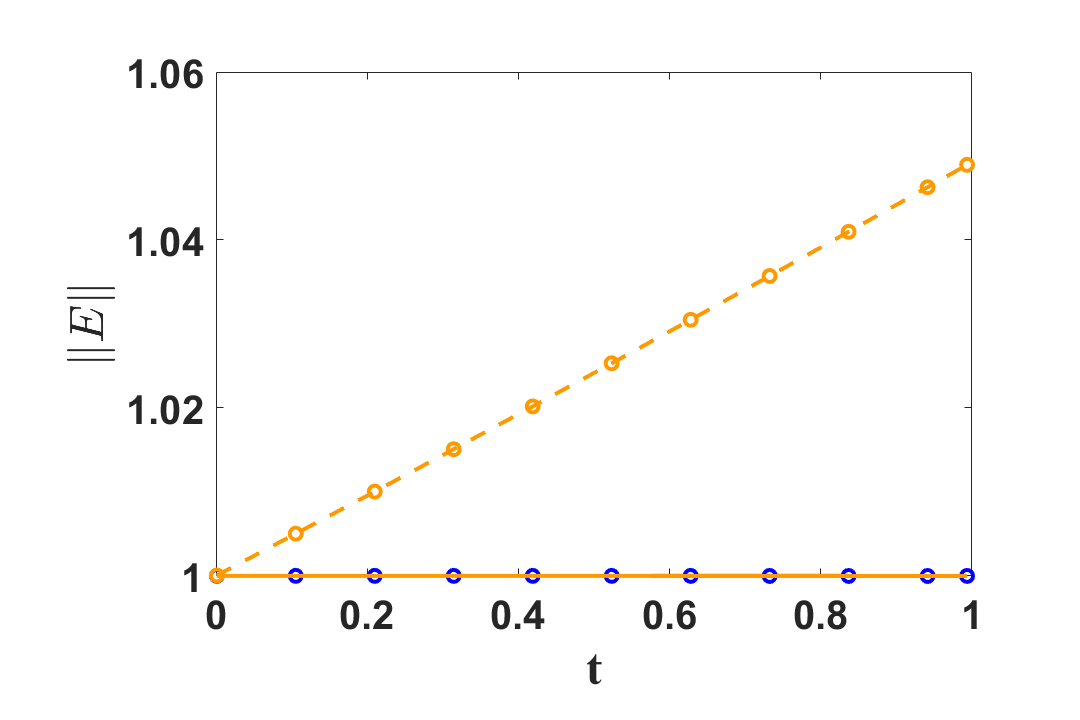}
} 
\mbox{
\hspace{0.4cm}
\makebox[0.45\textwidth]{(c)}
 \hspace{0.05\textwidth}
\makebox[0.45\textwidth]{(d)}
}

\caption{Results for a two-dimensional constant velocity advection with DSEM-SL using $H$=4 subdomains in $x$- and $y$-direction and an approximation order of $P$=6 per subdomain; (a) the solution, $\phi$, versus $x$  at $T$=1. 
The $L^2$ error norm, $\|e\|_{L^2}$, the mass norm, $\|M\|$, and the energy norm, $\|E\|$ are plotted versus time, $t$, in subfigures (b), (c) and (d), respectively}  
\label{fig:3_1}
\end{figure}

 The time evolution of the conservative properties 
 is shown in Figure \ref{fig:3_1}. Since the velocity is constant, the time integration is exact in Basecase and Basecase2 and their solutions overlap each other, like for
 the one-dimensional case. 
 
\begin{figure} 
\centering 
\mbox{ 
  \includegraphics[width=0.5\textwidth]{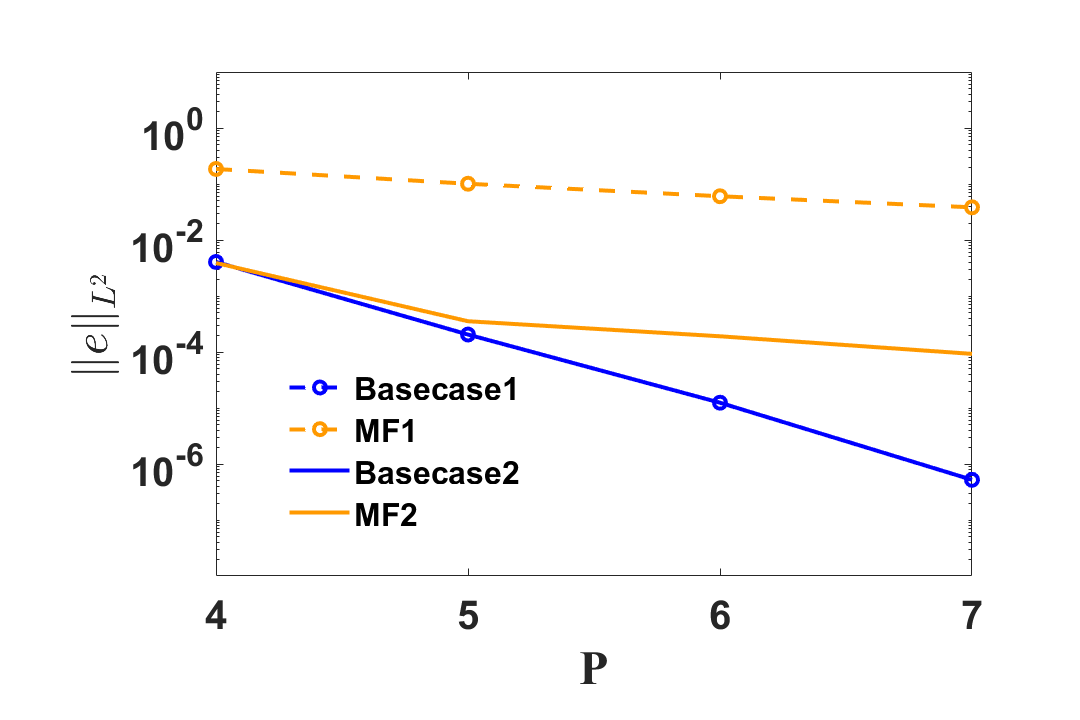} 
  \hspace{0.01\textwidth}
  \includegraphics[width=0.5\textwidth]{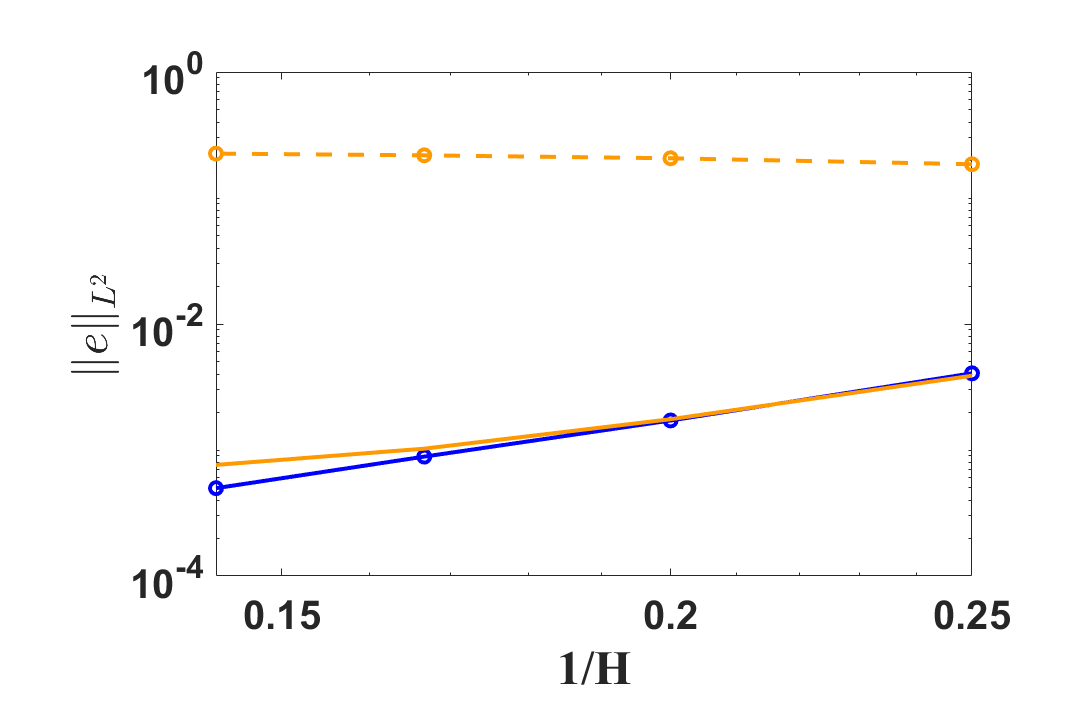}
}
\mbox{
\hspace{0.4cm}
\makebox[0.45\textwidth]{(a)}
 \hspace{0.05\textwidth}
\makebox[0.45\textwidth]{(b)}
}
\caption{$P$ convergence of the constant velocity 2D simulation using $H_x$=4, $H_y$=4 at $t$=1.} 
\label{fig:3_2}
\end{figure}
 Addition of the local mass constraints with a first order time integration MF1, which introduces an error of order $\mathcal{O}(10^{-3})$ to the system of equations (\ref{matrix_a}) leads to a global $L^2$ error of $\|e\|=\mathcal{O}(10^{-1})$. The MF1 scheme is non-conservative and gains $5\%$ energy over the time interval considered.
 
 The second order method MF2, with a time integration error on the order of $\mathcal{O}(10^{-6})$, is more accurate with the $L^2$ error on the order of $\|e\|=\mathcal{O}(10^{-4})$. 
 
 Figure \ref{fig:3_2} shows that the $P$ and $H$ convergence of the cases at $T$=1 are similar to the 1D results. Basecase1 and Basecase2 show spectral $P$-convergence as well as algebraic $H$ convergence.


 

\section{Conclusions}
\label{sec:5}
A high-order, explicit semi-Lagrangian method is developed
to solve the Lagrangian transport equations in Eulerian-Lagrangian formulations.
The semi-Lagrangian method is consistent with an explicit Eulerian solver
discretized with an explicit discontinuous spectral element method (DSEM).

By seeding tracer particles at the Gauss quadrature nodes, the Lagrangian
solution is directly available at quadrature nodes of the Eulerian solver and vice-versa. This is exchange of information is commonly performed using computationally intensive and complicated interpolation methods in Eulerian-Lagrangian formulations.

Consistent with DSEM, the semi-Lagrangian method is explicit. By
choosing the explicit time step appropriately,  particles are prevented from leaving the element. This ensures a
local and parallel method, which is natural for DSEM.

Following the explicit trace, the solution is remapped to the original
quadrature points using a least squares fit. This causes
the method to loose formal conservation.

Mass and kinetic energy constraints were tested to correct for this conservation loss.
The addition of a mass constraint  does not improve
mass conservation for cases with a constant advection velocity. For a non-constant advection velocity case,
however, the mass constraints can improve conservation. In general, an increase of
the time integration accuracy improves conservation and accuracy.

A kinetic energy constrains leads to an overconstraining of the system.
The solution shows spurious oscillations and is unstable.

We are extending the current algorithm to solve stochastic differential 
equations in Eulerian-Lagrangian formulations and plan to report on this in the near future.


\section*{Acknowledgements}
Funding provided by the Computational Science Research Center and AFOSR under grant number FA9550-16-1-0008 and NSF-DMS 1115705 is greatly appreciated.

\bibliographystyle{siam}
\bibliography{main}

\end{document}